\PassOptionsToPackage{dvipsnames}{xcolor}

\documentclass[lettersize,journal]{IEEEtran}
\usepackage{amsmath,amsfonts}
\usepackage{algorithmic}
\usepackage{algorithm, setspace}
\usepackage{array}
\usepackage[caption=false,font=normalsize,labelfont=sf,textfont=sf]{subfig}
\usepackage{textcomp}
\usepackage{fixltx2e}
\usepackage{dblfloatfix}
\usepackage{float}
\usepackage{url}
\usepackage{verbatim}
\usepackage{graphicx}
\usepackage{cite}
\usepackage{bbold}

\usepackage{caption}
\usepackage{textcomp}
\usepackage{xcolor}
\usepackage{amssymb}
\usepackage{mathtools}
\usepackage{color}
\usepackage[utf8]{inputenc}
\usepackage[english]{babel}
\usepackage{amsthm}
\newtheorem{theorem}{Theorem}[section]

\usepackage{comment}

\newcommand{\bx}{\mathbf{x}}
\newcommand{\bF}{\mathbf{F}}
\newcommand{\bG}{\mathbf{G}}
\newcommand{\bT}{\mathbf{T}}
\newcommand{\bI}{\mathbf{I}}

\newcommand{\grad}{\nabla}

\newcommand{\bv}{\mathbf{v}}

\DeclareMathOperator*{\argmin}{arg\,min}

\renewcommand{\comment}[1]{}

\newcommand{\extraproof}[1]{} 

\begin{document}

\title{Projected Multi-Agent Consensus Equilibrium (PMACE) with Application to Ptychography}

\author{Qiuchen Zhai,~\IEEEmembership{Student Member,~IEEE,} Gregery T. Buzzard,~\IEEEmembership{Senior Member,~IEEE,} Kevin Mertes, Brendt Wohlberg,~\IEEEmembership{Fellow,~IEEE,} Charles A. Bouman,~\IEEEmembership{Fellow,~IEEE}
\thanks{
The work of Q. Zhai was supported by U.S. Department of Energy through LANL.  
The work of C. A. Bouman was supported in part by the U.S. Department of Energy and the Showalter Trust. 
The work of G.T. Buzzard was partially supported by NSF CCF-1763896. 
The work of B. Wohlberg was supported by the LDRD program of LANL under project number 20200061DR. 

Qiuchen Zhai and Charles A. Bouman are with the School of Electrical and Computer Engineering, Purdue University, 465 Northwestern Ave., West Lafayette, IN 47907, USA (qzhai@purdue.edu and bouman@purdue.edu). 

Gregery T. Buzzard is with the Department of Mathematics, Purdue University, West Lafayette, IN 47907, USA  (buzzard@purdue.edu).

Kevin Mertes is with Physical Chemistry and Applied Spectroscopy Group, Los Alamos National Laboratory, Los Alamos, NM 87545 USA (kmmertes@lanl.gov).

Brendt Wohlberg is with Theoretical Division, Los Alamos National Laboratory, Los Alamos, NM 87545 USA (brendt@ieee.org).}}

\maketitle

\begin{abstract}
Multi-Agent Consensus Equilibrium (MACE) formulates an inverse imaging problem as a balance among multiple update agents such as data-fitting terms and denoisers. 
However, each such agent operates on a separate copy of the full image, leading to redundant memory use and slow convergence when each agent affects only a small subset of the full image.  

In this paper, we extend MACE to Projected Multi-Agent Consensus Equilibrium (PMACE), in which each agent updates only a projected component of the full image, thus greatly reducing memory use for some applications.  
We describe PMACE in terms of an equilibrium problem and an equivalent fixed point problem and show that in most cases the PMACE equilibrium is not the solution of an optimization problem. 

To demonstrate the value of PMACE, we apply it to the problem of ptychography, in which a sample is reconstructed from the diffraction patterns resulting from coherent X-ray illumination at multiple overlapping spots. 
In our PMACE formulation, each spot corresponds to a separate data-fitting agent, with the final solution found as an equilibrium among all the agents. 
Our results demonstrate that the PMACE reconstruction algorithm generates more accurate reconstructions at a lower computational cost than existing ptychography algorithms when the spots are sparsely sampled.

\end{abstract}


\section{Introduction}
\IEEEPARstart{M}{ethods} for inverse imaging typically seek to balance a fit to noisy data with some form of regularization or prior information.  
While the fit to data is often well-understood from physical models and can be captured using a log-likelihood or other penalty function, the prior information is much more difficult to capture in such a functional form.  
On the other hand, algorithmic denoisers such as BM3D~\cite{dabov2007} and convolutional neural network (CNN) denoisers~\cite{zhang2017beyond} encapsulate a great deal of prior information about images but only in an algorithmic way that is not easily amenable to description as a penalty function.  

The recently developed Multi-Agent Consensus Equilibrium (MACE)~\cite{buzzard2018plug} and Plug-and-Play (PnP) \cite{venkatakrishnan_plug-and-play_2013, sreehari2016plug, bouman2022foundations} methods provide a formulation and algorithms for incorporating both algorithmic prior information and information from penalty functions into a single reconstruction framework.  
Because the algorithmic prior information is encoded in input-output format as in an image denoiser, there is no cost function associated with this prior information, hence no cost function and no minimization associated with the reconstruction.  
Instead, MACE formulates the problem as an equilibrium among multiple input-output maps called agents.  

The theory underlying MACE implies that PnP methods are equivalent to Bayesian methods when the agents are proximal maps of convex functions (described below), but that PnP methods are more general in that there is no corresponding natural optimization problem when the PnP agents are not proximal maps, as in the case of neural networks. Moreover, MACE generalizes PnP to incorporate more than 2 agents, which allows the problem to be divided into smaller subproblems and solved separately with distributed implementation.  However, the original MACE framework requires that each agent maintain a separate copy of the entire reconstruction, which is very inefficient when each agent operates on only a small portion of the reconstruction.

Ptychography is an important computational imaging technique~\cite{rodenburg2008ptychography, rodenburg2019ptychography} that is widely used in applications such as the imaging of manufactured nanomaterials~\cite{wilke2012hard, trtik2013density, guizar2014high, shapiro2017ptychographic}.
The imaging technique works by moving a coherent X-ray probe across an object plane to scan a sample in overlapping patches.
A detector then measures the magnitude of the far-field diffraction pattern at each probe location. The magnitude and phase of the transmittance image of the sample can be recovered by exploiting the redundancy caused by the overlaps in the intensity measurements.  Notably, the reconstructed phase is significant in estimating the structure
of the sample.

Since ptychography is a phase recovery problem, it can be solved using conventional phase recovery methods~\cite{fienup1982phase}.
However, a series of methods have been developed specifically for reconstruction of ptychographic data. 
The original method, known as ptychographic iterative engine (PIE)~\cite{rodenburg2004phase}, has been generalized to methods such as extended PIE (ePIE)~\cite{maiden2009improved}, regularized PIE~(rPIE)~\cite{thibault_reconstructing_2013}, and momentum-accelerated PIE (mPIE)~\cite{maiden2017further}.
All these methods tend to have rapid convergence and work by serially enforcing local intensity measurements.
However, since the PIE-type algorithms are based on serial updates of the image at each probe location, the PIE methods cannot be easily parallelized, and the methods also tend to require a large overlap between scan positions to achieve a high-quality reconstruction~\cite{bunk2008influence}.

In addition, a variety of gradient-based methods including Wirtinger Flow (WF)~\cite{candes2015phase}, accelerated WF~\cite{xu2018accelerated}, Maximum Likelihood method~\cite{odstrvcil2018iterative}, as well as other approaches based on 
proximal algorithms~\cite{soulez2016proximity, yan2020ptychographic} have also been proposed for ptychographic reconstruction.
While methods based on gradient calculations can be parallelized across probe locations, these methods tend to have slower convergence than serial methods and are limited by the requirement that the problem must be formulated as an optimization.

More recently, a method known as scalable heterogeneous adaptive real-time ptychography (SHARP) has been proposed~\cite{marchesini2016alternating} and  implemented~\cite{marchesini2016sharp} for ptychographic reconstruction that is intrinsically parallel in nature.
The SHARP method is based on a general method for phase-retrieval known as relaxed averaged alternating reflections (RAAR)~\cite{luke2004relaxed}.
However, SHARP and RAAR are designed specifically for phase retrieval and not more general inverse problems.

In this paper, we introduce Projected Multi-agent Consensus Equilibrium (PMACE) an intrinsically parallel algorithm for solving distributed inverse problems.
PMACE extends the MACE framework to allow the state of a problem to be updated in smaller components that can be processed in parallel by localized agents. 
For some applications such as ptychography, PMACE can dramatically reduce memory requirements relative to MACE because each agent operates in parallel while storing only a small portion of the entire image.
We also show that PMACE is equivalent to an optimization formulation for some choices of update agents, but that in most cases there is no naturally corresponding optimization problem, with the solution corresponding instead to a weighted equilibrium of the updates of the individual agents. 

To illustrate the utility of the PMACE framework, we apply it to the problem of ptychographic reconstruction.   
We use the data from each scan position to implement a computationally efficient, non-optimization-based data-fitting agent that updates a local patch to better fit the measured intensity data for a single diffraction spot. 
These local patches are combined into a full image with a weighted average that uses the probe intensity to capture the uncertainty in measurement at a given location.  

In contrast to ePIE, the inherent parallelism of individual agents allows the PMACE algorithm for ptychography to make parallel updates to the image at different probe locations.  
Therefore, as with SHARP, the PMACE algorithm is more suitable than PIE-type algorithms for implementation on parallel or distributed computers.  However, we leave the exploration of efficient parallel implementation for future work.

Our experimental results demonstrate that PMACE outperforms existing state-of-the-art ptychographic reconstruction algorithms in both convergence speed and reconstruction quality as the ptychographic spots are spaced further apart. In practice, reconstruction from sparse samples is of practical importance because it can be used to reduce the amount of data required for accurate ptychographic reconstruction.

\IEEEpubidadjcol

\section{MACE and PMACE}

\subsection{MACE}

Multi-Agent Consensus Equilibrium (MACE)~\cite{buzzard2018plug} is a problem formulation that reconciles multiple agents, such as data fitting updates or denoisers, each acting to improve a candidate reconstruction.  
In this formulation, each agent, $F_j$, maintains an individual copy, $x_j$, of the full image (or volume), and $F_j(x_j)$ is an improved reconstruction according to the $j$th agent.  
These $x_j$ are either real or complex vectors with a dimension appropriate to the problem, and $F_j(x_j)$ has the same dimension.  We leave the dimensions implicit when they are clear from the context.
These copies are stacked to form a full MACE state, given by
$$
\bx = 
\left[
x_0,
\ldots,
x_{J-1} 
\right] \ ,
$$
and the agents are stacked to form a single operator defined by
\begin{equation}
\label{eq: stacked_forward_operators_x}
    \bF (\bx) = 
    \left[
    F_{0} (x_{0}), 
    \ldots,
    F_{J-1} (x_{J-1})
    \right] \ . 
\end{equation}

Note that $\bx$ denotes the full MACE state, while $x_j$ denotes the $j$th component image. 
Since each component is a full copy of the image to be reconstructed, the MACE state contains multiple, potentially inconsistent reconstructions.
In order to produce a single MACE reconstruction, we define the operator
\begin{equation}
\label{eq:consensus_operator_x}
    \bG ( \bx ) =
    \left[
    \bar{x}(\bx ),
    \ldots,
    \bar{x}(\bx )
    \right] \ , 
\end{equation}
where $\bar{x}(\bx ) = (1/J) \sum_{j=0}^{J-1} x_j.$
So $\bG ( \bx )$ computes the average of the components in $\bx$, and then returns a state vector formed by replicating this average $J$ times.
A key algebraic property of the averaging operator is that $\bG ( \bG (\bx ))= \bG (\bx )$. 

The MACE equation is then given by
\begin{equation} \label{eq:MACE-equil}
\bF(\bx^*) = \bG(\bx^*) \ , 
\end{equation}
where $\bx^*$ solves the equation, and the final reconstruction is then given by $x^* = \bar{x}(\bx^*)$.  As shown in~\cite{buzzard2018plug}, \eqref{eq:MACE-equil} has the interpretation of finding a consensus equilibrium among all the agents, in that the output of each agent is the common image $\bar{x}(\bx^*)$, while the updates $F_j(x_j^*) - x_j^*$ sum to 0.  

In the case that each $F_j$ is the proximal map for a convex, differentiable\footnote{We assume differentiability for simplicity.  Similar results hold for general convex functions by using subgradients.} function $f_j$, then $x^*$ satisfies 
\begin{equation} \label{eq:MACE-grad-zero}
   \sum_{j=0}^{J-1} \grad f_j(x^*) = 0 \ ,
\end{equation}
and $x^*$ is a minimizer of the function $f = \sum_j  f_j$~\cite{sridhar2020distributed}.  

However, when the MACE agents are not proximal maps, there may be no cost function that is minimized by solving the MACE equation. 
We discuss this in greater generality in Section~\ref{sec:PMACE}.

We note that $\bF (\bx )$ is an intrinsically parallel operation that can be efficiently distributed across many compute nodes, while the averaging operator $\bG ( \bx )$ requires communication across all nodes in order to first gather the individual components, then compute their average, and then broadcast the average back to individual nodes~\cite{sridhar2020distributed}.

\begin{figure}[t!]
    \centering
    \includegraphics[width=9cm]{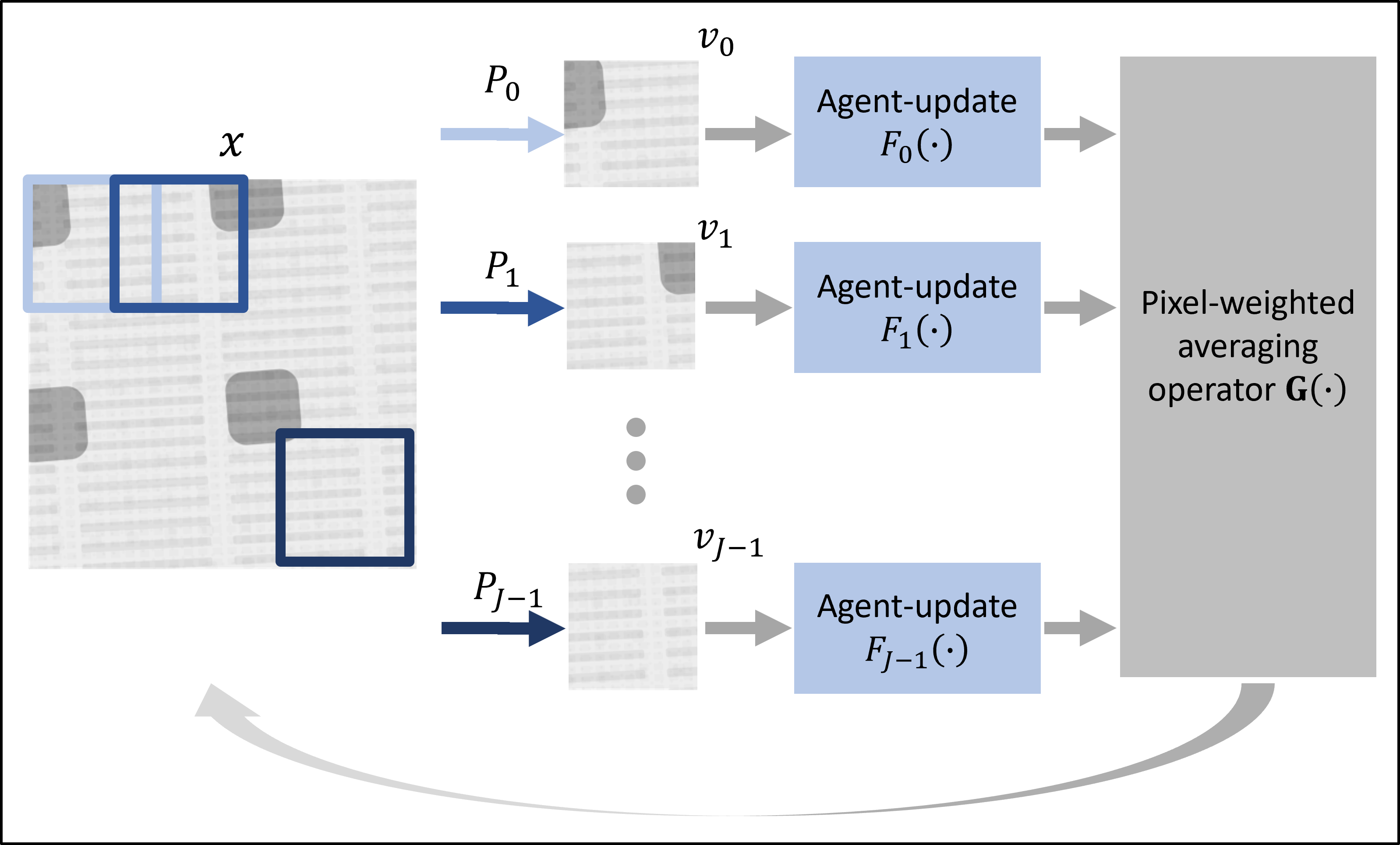}
    \caption{Conceptual rendering of the PMACE pipeline. A reconstruction $x$ is subdivided into possibly overlapping components, $v_j$, and distributed to multiple agents, $F_j$, each of which acts to improve this local reconstruction. These local reconstructions are combined using a pixel-weighted averaging operator, $G$, that reassembles the local components into a consistent global reconstruction.}
    \vspace{-0.6cm}
    \label{fig:pmace_flow}
\end{figure}

\subsection{PMACE}
\label{sec:PMACE}

In some problems, agents operate naturally on subsets or projections of a full image rather than on the image itself.
Figure~\ref{fig:pmace_flow} illustrates such a problem in which separate data is collected from each patch of an image, $x$.
In this case, the agents can operate separately to make each component more consistent with its data.
However, these individual components must also be reconciled into a single consistent reconstruction.

Mathematically, we express this extraction of a patch or component from an image as
$$
v_j = P_j x \ ,
$$
where $P_j$ is a matrix with orthonormal rows, so that $P_j P_j^T = I$.  
We use $v_j$ to emphasize that these components lie in a different space than $x$.  
Typically, $P_j$ simply selects a region of the image $x$, but more general transformations are allowed. To capture spatial information and relevant features without excessive computational cost, the selection of $P_j$ should be based on the size of local measurements associated with individual components.

As in MACE, each component, $v_j$, is processed by a corresponding agent, $F_j$, whose role is to improve the component by better fitting to data, reducing noise, etc.
Also as in MACE, PMACE stacks the component states and the agents to form the vector and operator given by
\begin{equation}
\label{eq: stacked_vj}
\bv = 
\left[
v_0,
\ldots,
v_{J-1} 
\right] \ ,
\end{equation}
\begin{equation}
\label{eq: stacked_forward_operators_v}
    \bF (\bv) = 
    \left[
    F_{0} (v_{0}), 
    \ldots,
    F_{J-1} (v_{J-1})
    \right] \ . 
\end{equation}

Note that PMACE has a potentially huge advantage over MACE because it does not require $J$ replicates of the entire reconstruction.
This can dramatically reduce memory and computation resources when implemented on a computer.
For example, in ptychography as described below, if each reconstruction is a large image and there are $1000$ agents corresponding to $1000$ probe locations, then MACE becomes impractical to implement since it requires that the full image be replicated $1000$ times in the computer's memory. 

Since the PMACE components are not complete versions of the reconstruction, the averaging operator $\bG (\bx )$ must be redefined.
To do this, we project the components back to the full image, then use a pixel-wise weighted average to combine them and reproject them back to components.  This leads to the new averaging operator
\begin{equation}
\label{eq:consensus_operator_v}
    \bG^P ( \bv ) =
    \left[
    P_0 \bar{x}(\bv),
    \ldots,
    P_{J-1} \bar{x}(\bv)
    \right] \ . 
\end{equation}
where 
\begin{equation}
\label{eq:consensus_operator_component_v}
\bar{x}(\bv) = \Lambda^{-1}  \sum_{j=0}^{J-1} P_j^T W v_j \ ,
\end{equation}
${\Lambda} = \sum_{j=0}^{J-1} P_{j}^T W P_j $, and $W$ is a positive definite and (typically) diagonal weighting matrix.
We assume that $\Lambda$ is invertible, so if $P_j$ selects subsets of the original image, then the set of $J$ components must cover the entire image being reconstructed.
The weight matrix, $W$, can be used to model the uncertainty in measurements for each component, with its selection depending on the nature of imaging modality and the noise characteristics in the data.

Using these operators, we look for a solution $\bv^*$ to the PMACE equation, defined as 
\begin{equation} \label{eq:pmace-eqns}
    \bF(\bv^*) = \bG^P(\bv^*) \ ,
\end{equation}
in which case $x^* = \bar{x}(\bv^*)$ is the final reconstruction. In words, \eqref{eq:pmace-eqns} means that each agent generates a patch that is consistent with the entire reconstruction $x^*$.

\subsection{Solving the PMACE Equations}

To solve the PMACE equation \eqref{eq:pmace-eqns}, we follow~\cite{buzzard2018plug} to show that it can be reformulated as a fixed point problem and apply Mann iterations to compute a solution.

By definition of $\bar{x}$, $G^P$, and $\Lambda$, we have
\begin{align*}
    \bar{x}(\bG^P(\bv)) &= \Lambda^{-1} \sum_j P_j^T W P_j \bar{x}(\bv)\\
    &= \bar{x}(\bv) \ ,
\end{align*}
which gives 
\begin{equation}  \label{eq:GG=G}
    \bG^P ( \bG^P (\bv ))= \bG^P (\bv ) \ .
\end{equation}
Thus $\bG^P$, like the original averaging operator $\bG$, is a linear projection.
With this, a further calculation shows that 
\begin{equation} \label{eq:2G-I-inv}
(2 \bG^P -\bI)^{-1} = (2 \bG^P -\bI) \  .  
\end{equation}
From this we obtain the equivalent fixed point formulation of \eqref{eq:pmace-eqns} as 
\begin{equation}  \label{eq:T-fixed}
    (2 \bG^P -\bI)(2 \bF - \bI)(\bv^*)  = \bv^* \ ,
\end{equation}
or $\bT \bv^* = \bv^*$, where $\bT = (2 \bG^P -\bI)(2 \bF - \bI)$.

If $\bT$ is non-expansive and has a fixed point, then Mann iterations with $\rho \in (0,1)$, defined as  
\begin{equation}  \label{eq:Mann-T}
    \bv \gets (1 - \rho) \bv + \rho \bT \bv \ ,
\end{equation}
are guaranteed to converge to a fixed point of $\bT$ and hence to a solution of \eqref{eq:pmace-eqns}. The parameter $\rho$ is a kind of step size; its choice impacts the convergence rate and the stability of the iterations but not the final solution.

\subsection{PMACE for Optimization and Beyond} %
\label{subsec:PMACEOptimization}

A natural question is whether PMACE corresponds to the solution of an optimization problem.  We show in Theorem~\ref{thm: using_prox_equil_problem} that the answer is affirmative in certain special cases, but not in general. (Proofs are given in the appendix.)

To understand the relationship between equilibrium and optimization, suppose that each $F_j$ is the proximal map for a real-valued, differentiable, convex function $f_j(v_j)$:
\begin{equation}
\label{eq: proximal_map_function}
    \begin{aligned}
        F_{j} (v_{j}) = \arg \min _ {v} \left \{ f_{j}(v) + \frac{1}{2 \sigma^2} \left \| v - v_{j} \right \|^2 \right \} \ .
    \end{aligned}
\end{equation}

In this setting, a PMACE solution is equivalent to a zero of a naturally associated vector field.  
\begin{theorem}
\label{thm: using_prox_equal_optimization_problem}
Let $F_j, j = 0, \cdots, J-1$ denote the proximal map function of a differentiable convex function $f_j$ as specified in \eqref{eq: proximal_map_function}.  Then $\bv^*$ is a solution to the PMACE equation in \eqref{eq:pmace-eqns} if and only if $x^* = \bar{x}(\bv^*)$ satisfies
\begin{equation}
    \label{eq:vec-field}
    \sum_{j=0}^{J-1} P_j^T W \nabla f_j(P_j x^*) = 0 \ .
\end{equation}
\end{theorem}

Note that \eqref{eq:vec-field} captures the idea of equilibrium in that the sum of the $J$ vectors is 0.  
However, as described in more detail below, because of the weight matrix $W$, this vector field may not be a gradient field, hence may not correspond naturally to a minimization problem.  
To explain this more fully, we first recall some basic results on vector fields.  

We say that $V(x)$ is a conservative vector field if it is the gradient of some potential function $f(x)$~\cite[Section 9.2A]{williamson2004multivariable}.
In this case, the condition $V(x^*)=0$ is a first-order necessary condition for the minimization of $f$.
Moreover, a continuously differentiable vector field, $V(x)$, is conservative if and only if its Jacobian is self-adjoint (or symmetric)~\cite[Section 9.2C]{williamson2004multivariable}.  

When each $f_j$ is twice continuously differentiable, the Jacobian of the vector field in \eqref{eq:vec-field} (as a function of $x$ in place of $x^*$) is
\begin{equation}
    \label{eq:grad-vec-field}
    \sum_{j=0}^{J-1} P_j^T W H_{f_j}(P_j x) P_j \ ,
\end{equation}
where $H_{f_j}(P_j x)$ is the symmetric Hessian matrix of 2nd order partial derivatives of $f_j$, evaluated at $P_j x$.  
If $W$ is a multiple of the identity, then the matrix in \eqref{eq:grad-vec-field} is symmetric, so the PMACE solution corresponds naturally to a minimization problem.  
However, this symmetry is not robust in that nonzero mixed partial derivatives in $f_j$ together with distinct entries in the diagonal matrix $W$ will destroy this symmetry.  
This leads to the following theorem.
\begin{theorem}
\label{thm: using_prox_equil_problem}
Let $f_j$ and $F_j$ be as in Theorem~\ref{thm: using_prox_equal_optimization_problem}.  
\begin{enumerate}
    \item If $W = rI$ for scalar $r$,  then PMACE is equivalent to optimization in that $\bv^*$ is a solution to the PMACE equation in \eqref{eq:pmace-eqns} if and only if $x^* = \bar{x}(\bv^*)$ satisfies
    \begin{equation}
    \label{eq: ml_est}
    x^* = \arg\min_{x}\left \{\sum_{j=0}^{J-1}f_{j}(P_j x) \right \} \ .
    \end{equation}
    \item For generic diagonal $W$ and generic convex $f_j$, the PMACE formulation does not naturally arise as an optimization problem.  That is, the vector field in \eqref{eq:vec-field} is not a conservative vector field, and hence is not the gradient field of a potential function. 
\end{enumerate}
\end{theorem}

In other words, even when the agents $F_j$ are proximal maps, if the weight matrix $W$ is not a multiple of the identity, then the natural vector field defining the PMACE solution is not a gradient field, so the PMACE solution is not naturally given by minimizing a cost function.  
Of course, it's possible to set up an optimization problem with the same solution by minimizing the norm of the vector field in \eqref{eq:vector-field-V}, but this doesn't correspond to a consensus minimization problem as in \eqref{eq: ml_est}.  

When the $f_j$ are convex but not everywhere differentiable, then the gradient $\nabla f_j$ in \eqref{eq:vec-field} is replaced with the possibly set-valued subdifferential $\partial f_j$, and the search for an $x^*$ to give a zero of the left-hand side becomes the search for an $x^*$ so that the left-hand side contains the 0 vector, which is known as a monotone inclusion problem. 
Results analogous to Theorem~\ref{thm: using_prox_equil_problem} can be obtained in this setting but are not pursued here.  
For more about the relationship between convex optimization and monotone operators, see~\cite{combettes_monotone_2018}.     

\bigskip
{\bf Non-proximal equilibrium:}  
To extend to the case in which the agents $F_j$ are not proximal maps, we note that the proximal map in \eqref{eq: proximal_map_function} can be rewritten~\cite[Cor.~17.6]{Bauschke2011ConvexAA} as $F_j(v_j) = v_j - \sigma^2 \nabla f_j(F_j(v_j))$.  
That is, a proximal map for a convex, differentiable function is an implicit gradient descent step in the sense that the gradient is evaluated at the endpoint rather than the starting point.  
With this and the fact that $F_j(v_j^*) = P_j x^*$, \eqref{eq:vec-field} can be rewritten as 
\begin{equation} \label{eq:PMACE-optimization-equilibrium-v}
     \sum_{j=0}^{J-1} P_j^T W (F_j(v_j^*) - v_j^*) = 0 \ .
\end{equation}

The vector field in this case is defined not on an image $x$ but on a stack of patches as in \eqref{eq: stacked_vj}.  In order to use this formulation to give a single reconstruction, $x^*$, we note that each $v_j^*$ must satisfy the constraint $F_j(v_j^*) = P_j x^*$.  We use this relationship to implicitly define a set-valued inverse $F_j^{-1}(P_j x))$ that maps $x$ to all patches $v_j$ consistent with $x$.   

With this, we reformulate \eqref{eq:PMACE-optimization-equilibrium-v} and hence the PMACE equation in \eqref{eq:pmace-eqns} as an inclusion problem.

\begin{theorem} 
\label{thm: pmace_sol_equal_v_goes_zeros}
Define the set-valued vector field 
\begin{equation}  \label{eq:vector-field-V}
   V(x) = \sum_{j=0}^{J-1} P_j^T W \left( P_j x  - F_j^{-1}(P_j x)\right) \ . 
\end{equation}
Then the solution $x^* = \bar{x}(\bv^*)$ of the PMACE equation in ~\eqref{eq:pmace-eqns} is the solution of the inclusion problem 
\begin{equation}  \label{eq:inclusion-V}
0 \in V(x^*) \ .
\end{equation}
\end{theorem}

We call this a consensus equilibrium problem because the common point $x^*$ yields a set of vectors $v_j^*$ with $F_j(v_j^*) = P_j x^*$ (consensus) such that the weighted average updates give the zero vector as in \eqref{eq:PMACE-optimization-equilibrium-v}.  

The verification that a given operator is monotone is often technically involved, so we do not address that here. Likewise, further conditions under which \eqref{eq:vector-field-V} corresponds to an optimization problem are beyond the scope of this paper.   

Finally, we note that, as with MACE or PnP \cite{venkatakrishnan2013plug}, PMACE can incorporate regularization by either using a denoiser as one of the agents, or using a denoiser as part of the averaging operator, $\bG (\bx )$.
While this could potentially improve results for ptychography reconstruction, we leave this topic as an area for future research.

\section{PMACE for Ptychography}
\label{sec:PMACEPtychography}

\begin{figure}[htb]
    \centering
    \vspace{0.4cm}
    \includegraphics[width=9cm]{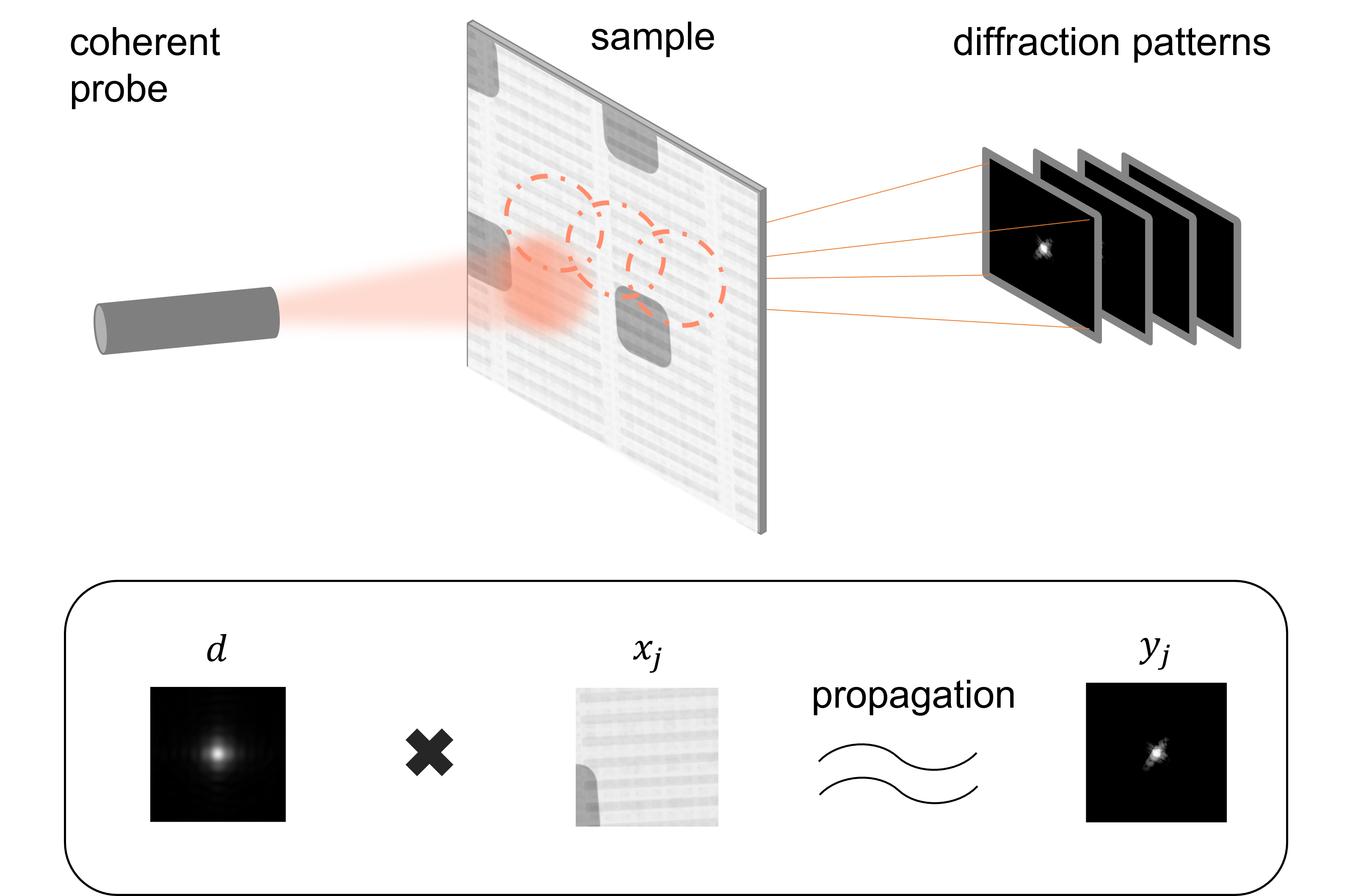}
    \vspace{0.1cm}
    \caption{In ptychographic imaging, a coherent probe is used to illuminate a sample at multiple, partially overlapping regions. The intensities of the resulting diffraction patterns are recorded by the detector.}
    \label{fig:ptycho_illustration}
\end{figure}

As noted in the introduction and illustrated in Figure~\ref{fig:ptycho_illustration}, a ptychographic reconstruction relies on far-field intensity measurements (also known as diffraction patterns) obtained by illuminating a sample with a coherent X-ray source at multiple locations, with significant overlap between adjacent illuminated regions.    
The goal of these measurements is to reconstruct the full complex transmittance image of the scanned sample.  
The phase of this complex transmittance image is particularly important since it gives precise information about the thickness of the sample at each location.  

A potential point of confusion in this problem is that there are two sets of phases that must be estimated for this reconstruction.  One set is the phase of the complex transmittance image, which is the primary quantity of interest because it gives very precise, quantitative information about the material properties of the sample.  The other set contains the phases associated with the far-field intensity data.  These far-field phases are essentially the phases of the Fourier transform of the illuminated sample.  Since the detector measures intensities only, the far-field phases must be estimated using the redundancies associated with overlapping intensity measurements -- this is known as the phase-retrieval problem.  However, for our purposes, the far-field phases are really just nuisance parameters that must be estimated in order to estimate the phase and magnitude of the complex transmittance.  

To apply PMACE to the ptychography problem, we define each data-fitting agent to be a local update to better match the measurements of the far-field Fraunhofer diffraction pattern in a single patch.  
Our approach does not use a penalty function/proximal map formulation but rather defines an update agent for each patch update directly from data using a weighted average between a previous estimate and new estimate that fits the available data.  This approach has the advantages of being computationally efficient and easily adjusted to account for the signal-to-noise ratio of the data.      

While probe estimation is an important component of ptychography, we assume for clarity in this paper that an accurate estimate of the probe function is known and consistent at all scanning positions.  
In a future paper we will apply PMACE to estimate both the probe and complex image simultaneously.

\subsection{Forward Model}

We let $x \in \mathbb{C}^{N_{1} \times N_{2}}$ denote the complex material transmittance of the target object and seek to recover $x$ from real-valued measurements $y_{j} \in \mathbb{R}^{N_{p}\times N_{p}}$ of the form
\begin{equation}
\label{eq: mag_measurement}
y_{j} = \sqrt{\text{Pois} \left ( | \mathcal{F} D P_{j} x |^2 \right )} \ .
\end{equation}
Here $\mathcal{F}$ denotes the 2D orthonormal discrete Fourier transform, $D = \text{Diag}(d)$ is a diagonal matrix representing the probe's complex illumination function $d$, and $P_{j}: \mathbb{C}^{N_{1}\times N_{2}} \rightarrow \mathbb{C}^{N_{p}\times N_{p}}$ is the linear operator representing the extraction of an $N_{p}\times N_{p}$ patch corresponding to the $j^{th}$ scan position, $j=0,\ldots, J-1$.
Thus $\mathcal{F} D P_{j} x$ is the mean complex far-field diffraction pattern associated with the $j^{th}$ scan position.  
We assume that the physical sensor measurements can be modeled by a vector of independent Poisson distributed random variables denoted by $\text{Pois} \left ( \lambda \right )$ with mean $\lambda$ given by the pixel-wise intensity of this diffraction pattern, and that these measurements are pre-processed by taking a square-root, which serves as a variance-stabilizing transform for the Poisson distribution~\cite{godard2012noise}.

\subsection{PMACE Formulation}

In order to use the PMACE algorithm, we will need to select a data-fitting agent, $F_{j}(v_j)$. 
The proximal map is not a good choice because its evaluation requires the solution of a non-convex optimization problem, which is computationally expensive and not guaranteed to achieve a global minimum.
Fortunately, the PMACE method is not restricted to the use of proximal maps.
So instead, our strategy will be to design a data-fitting agent that moves the solution in a direction that better fits the data.

Using this approach, we define the data-fitting agents $F_{j}(v_j)$ by taking a weighted average between a nonlinear projection onto measured data and an input point.  
The nonlinear projection begins with a complex-valued image patch, $v_j= P_j x$, and isolates the phase in diffraction space by $\mathcal{F} D v_j / |\mathcal{F} D v_j|$ (the phase can be set arbitrarily if the magnitude is 0).  This phase is multiplied point-wise by the amplitudes $y_j$, and then the result is inverted to return to image space.  

Mathematically, we represent the full nonlinear projection from an input image patch to one that matches the given amplitudes in diffraction space is given by $D_\epsilon^{-1}\mathcal{F}^{*} \left ( y_{j} \frac{\mathcal{F} D v_{j}}{|\mathcal{F} D v_{j}|} \right )$ where $D_\epsilon^{-1}=\mbox{Diag} (d_\epsilon^{-1})$ is a numerically stable inverse of $D$.\footnote{We compute the numerically stable inverse of the diagonal entries as $d_\epsilon^{-1} = d^* / (|d|^2+\epsilon)$ where $\epsilon=10^{-6} \sqrt{ \|d\|^2 /\mbox{dim}(d)} $.}

For the full data-fitting agent, we average this projection with the input point to obtain 
\begin{align} \label{eq:Fjv}
    F_{j} (v_{j})
        & = (1 - \alpha ) v_{j} + \alpha D_{\epsilon}^{-1} \mathcal{F}^{*} \left ( y_{j} \frac{\mathcal{F} D v_{j}}{|\mathcal{F} D v_{j}|} \right ) \ ,
\end{align}
where $\alpha$ controls the strength of fit to data. In contrast with $\rho$ in \eqref{eq:Mann-T},  where $\rho$ is employed to optimize convergence speed, $\alpha$ regulates regularization, with the ratio $\alpha / (1-\alpha)$ provides a rough interpretation of the signal-to-noise ratio of the data.

This interpolation between input and a data-fitting point is similar in spirit to a data-fitting proximal map. 
However, the Jacobian of \eqref{eq:Fjv} is not symmetric, so $F_j$ is not a proximal map (where we convert complex images to real-valued vectors by stacking real and imaginary parts in order to take real derivatives).  
In fact, Theorem~\ref{thm: pmace_sol_equal_v_goes_zeros} implies that the PMACE formulation using these $F_j$ does not naturally correspond to a minimization problem.  In the appendix, we show that $F_j$ is invertible, which implies that $V(x)$ in \eqref{eq:vector-field-V} is a single-valued vector field, which is nonconservative since the Jacobian of \eqref{eq:Fjv} is not symmetric.  Hence, as in Theorem~\ref{thm: using_prox_equil_problem}, this formulation is not naturally a minimization problem.

We use these $F_j$ to define the vectorized operator $\bF$ as in \eqref{eq: stacked_forward_operators_v}.  For $\bG$ in \eqref{eq:consensus_operator_v}, we  model the uncertainty arising from the point-to-point variation of probe illumination by using the weight matrix 
\begin{equation} \label{eq:Wkappa}
    W = |D|^\kappa
\end{equation}
to define $\Lambda$ and $\bar{x}$.  Here $\kappa \in [1,2]$ can be used to control the emphasis between brightly lit points near the center of the probe (large $\kappa$) versus dimly lit points near the edges of the probe (small $\kappa$).  
In our experiments, $\kappa \in [1.25, 1.5]$ gives the fastest convergence and best reconstruction quality; we display all results using $\kappa = 1.5$.  

The PMACE solution is then determined by solving $\bF(\bx) = \bG(\bx)$ or by iterating \eqref{eq:Mann-T} to obtain a fixed point of $\bT$. This has the effect of balancing the updates of the data fitting agents for each patch, subject to constructing an image that agrees with each individual patch.

Proper initialization of the PMACE algorithm is important for fast convergence to a good solution.
This is particularly true since phase retrieval corresponds to a non-convex optimization problem, so that solutions can in general depend on the initial condition.
We have found the following initialization to work well in practice, 
\begin{equation}
\label{eq:Initialization}
    x^{(0)} \gets \Lambda_{0}^{-1} \sum_{j=0}^{J-1} P_j^T \left(  \frac{ \| y_j \| }{ \| d \| } \mathbb{1} \right) \ ,
\end{equation}
where $\Lambda_{0} = \sum_{j=0}^{J-1} P_{j}^T P_j $ and $\mathbb{1}$ is a column vector of 1s. 
While simple, this creates an initial image with the correct scale at each probe location based on the relative magnitude of the probe and measured signal.

The resulting algorithm has similarities to SHARP~\cite{marchesini2016sharp} but also significant differences. 
First, equation \eqref{eq:Fjv} incorporates the data projector of~\cite{marchesini2016sharp}, but \eqref{eq:Fjv} is in the image domain and averaged with the identity, as opposed to the SHARP data projector, which is in the so-called frame domain (probe times image) and which is not averaged.  
The averaging factor $\alpha$ replaces the stabilization factor $\beta$ in the data projector of SHARP.  
More significantly, the averaging operator $P_Q$ of SHARP uses a fixed weighting that is essentially equivalent to taking $\kappa = 2$ in \eqref{eq:Wkappa}, which tends to underweight the contributions of locations with small to moderate probe intensities.

The pseudo-code for computing the PMACE solution is shown in Algorithm \ref{alg: PMACE}. 
The algorithm starts with an initial guess of the complex transmittance of the unknown object $x^{(0)} \in \mathbb C^{N_{1} \times N_{2}}$, and then constructs the stacked projections. Mann iterations are then used to revise the estimates of individual projections. The final reconstruction of the complex object is given by taking the weighted average of the estimates of individual projections.

\begin{algorithm}[H]
\setstretch{1.25}
\caption{Mann iteration for computing PMACE solution.}\label{alg: PMACE}
\begin{algorithmic}
 \renewcommand{\algorithmicrequire}{\textbf{Input:}}
 \renewcommand{\algorithmicensure}{\textbf{Output:}}
 \REQUIRE Initialization: $x^{(0)} \in \mathbb C^{{N_1} \times N_{2}}$; $\kappa \in [1.25, 1.5]$ 
 \ENSURE  Final Reconstruction: $\hat{x} \in \mathbb C^{{N_1} \times N_{2}}$
  \STATE $\mathbf{w} = \mathbf{v} = [v^{(0)}_{0},\dots, v^{(0)}_{J-1}],\; \text{where} \; v^{(0)}_{j} = P_{j} {x}^{(0)} $ \\
  \WHILE {not converged}
  \STATE \hspace{0.5cm}$\mathbf{w} \gets \bF ( \mathbf{v} )$
  \STATE \hspace{0.5cm}$\mathbf{z} \gets \bG (2\mathbf{w} - \mathbf{v})$
  \STATE \hspace{0.5cm}$\mathbf{v} \gets \mathbf{v} + 2 \rho (\mathbf{z} - \mathbf{w})$
  \ENDWHILE
 \RETURN $\hat{x} = \Lambda^{-1} \sum_{j=0}^{J-1} P_{j}^{T} |D|^{\kappa} v_{j}$ 
\end{algorithmic}
\end{algorithm}

\section{Experimental Results}
\label{sec: experimental result}
In this section, we present the results of our approach on both synthetic and measured data and compare with reconstruction results using state-of-the-art algorithms.  
Our results indicate that PMACE matches the best-competing algorithms in convergence speed and image quality for small probe spacing, and outperforms the competitors as the probe spacing between neighboring scan points increases. \footnote{ The code for PMACE is available at \url{https://github.com/cabouman/ptycho_pmace}.}

\subsection{Synthetic Data Generation}

\begin{figure}[!htp]
\vspace{-0.2cm}
\begin{minipage}[htb]{0.45\linewidth}
  \centering
  \hspace{0.1cm}
  \centerline{\includegraphics[height=3.4cm]{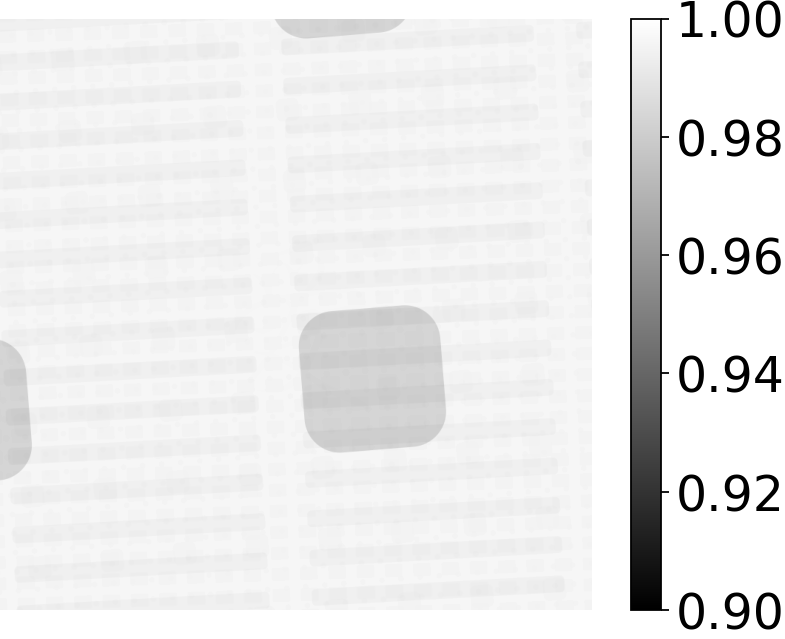}}
  \vspace{-0.2cm}
  \centerline{(a)} \smallskip
\end{minipage}
\begin{minipage}[htb]{0.45\linewidth}
  \centering
  \hspace{0.1cm}
  \centerline{\includegraphics[height=3.4cm]{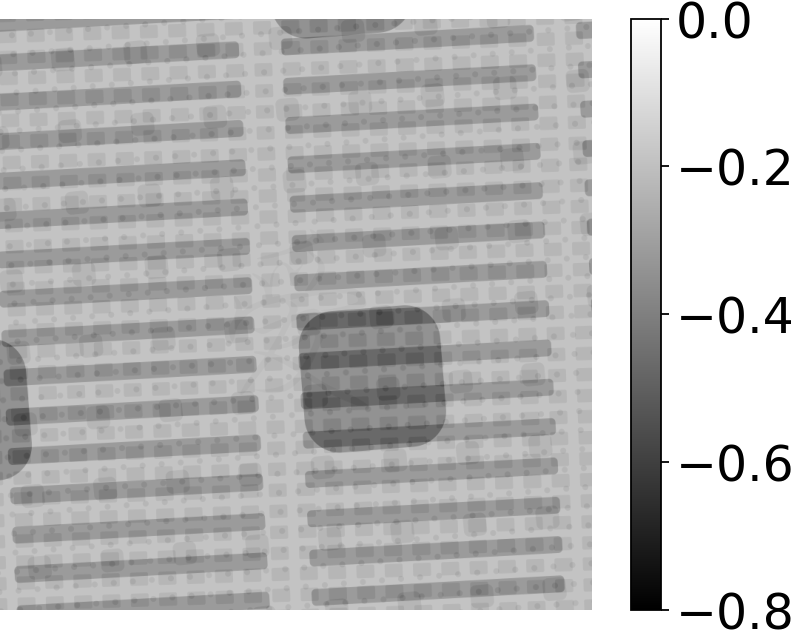}}
  \vspace{-0.2cm}
  \centerline{(b)} \smallskip
\end{minipage}
\hfill
\begin{minipage}[htb]{0.45\linewidth}
  \centering
  \hspace{0.1cm}
  \centerline{\includegraphics[height=2.7cm]{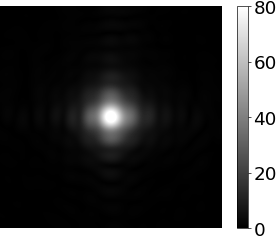}}
  \vspace{-0.2cm}
  \centerline{(c)} \smallskip
\end{minipage}
\hfill
\hspace{.6cm}
\begin{minipage}[htb]{0.45\linewidth}
  \centering
  \hspace{0.1cm}
  \centerline{\includegraphics[height=2.7cm]{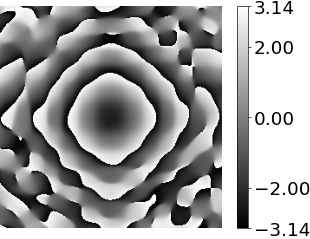}}
  \vspace{-0.2cm}
  \centerline{(d)} \smallskip
\end{minipage}
\hfill
\vspace{-0.1cm}
\caption{
Ground truth images used in simulation experiments.
The complex ground-truth object's (a) magnitude and (b) phase in radians.
The complex probe's (c) magnitude and (d) phase in radians.
}
\label{fig:ground_truth_imgs}
\end{figure}

Figure~\ref{fig:ground_truth_imgs} shows the complex-valued ground truth object image and probe used for our simulation experiments.
A synthetic complex transmittance image of size $800 \times 800$ pixels was generated by modeling a 5-layer 
composite material.

The complex probe of size $256 \times 256$ pixels was simulated with a photon energy of 8.8 keV.

\begin{figure*}[!htb]
\begin{minipage}[b]{.2\linewidth}
  \centering
  \centerline{\includegraphics[height=4cm]{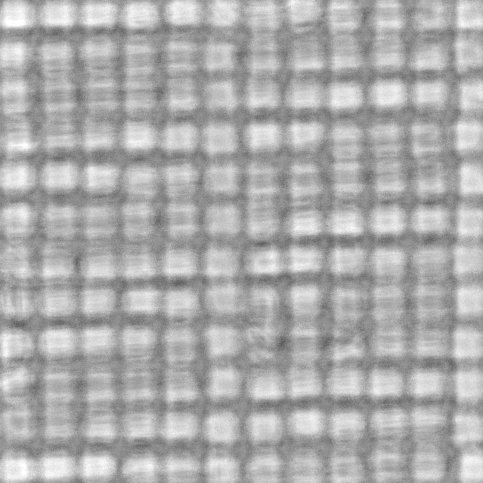}}
  \vspace{0.1cm}
  \centerline{(a) ePIE}\medskip
\end{minipage}
\hfill
\begin{minipage}[b]{.2\linewidth}
  \centering
  \centerline{\includegraphics[height=4cm]{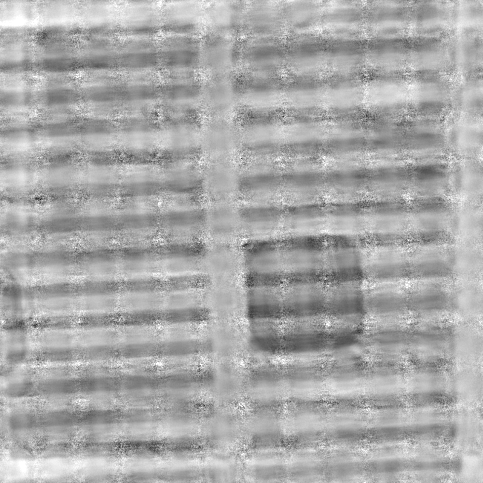}}
  \vspace{0.1cm}
  \centerline{(b) AWF}\medskip
\end{minipage}
\hfill
\begin{minipage}[b]{.2\linewidth}
  \centering
  \centerline{\includegraphics[height=4cm]{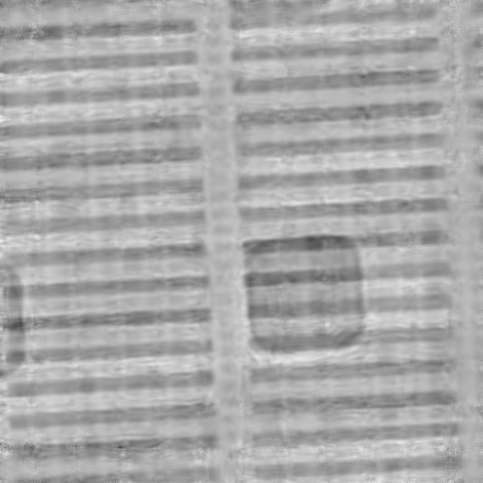}}
  \vspace{0.1cm}
  \centerline{(c) SHARP}\medskip  
\end{minipage}
\hfill
\begin{minipage}[b]{0.2\linewidth}
  \centering
  \centerline{\includegraphics[height=4cm]{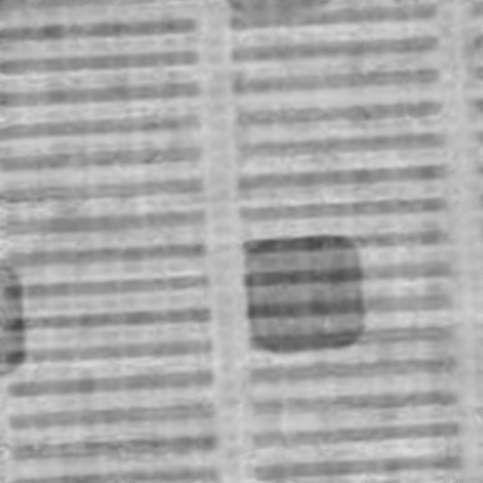}}
  \vspace{0.1cm}
  \centerline{(d) PMACE}\medskip
\end{minipage}
\hfill
\begin{minipage}[b]{.05\linewidth}
  \centering
  \centerline{\includegraphics[height=4cm]{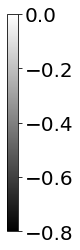}}
  \vspace{0.1cm}
  \centerline{ }\medskip  
\end{minipage}
\hfill
\begin{minipage}[b]{.2\linewidth}
  \centering
  \centerline{\includegraphics[height=4cm]{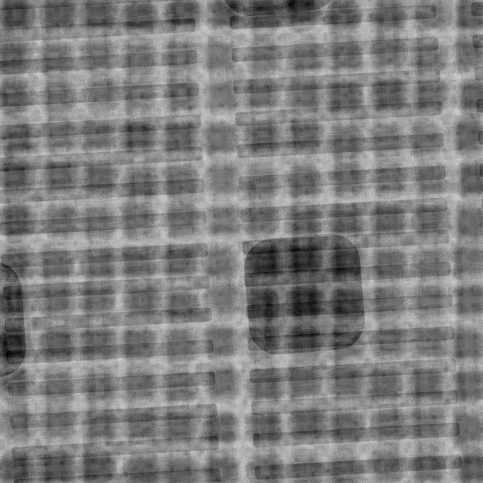}}
  \vspace{0.1cm}
  \centerline{(a) NRMSE = 0.175}\medskip
\end{minipage}
\hfill
\begin{minipage}[b]{.2\linewidth}
  \centering
  \centerline{\includegraphics[height=4cm]{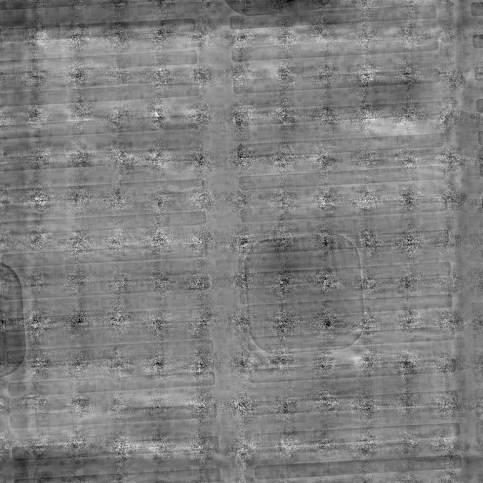}}
  \vspace{0.1cm}
  \centerline{(b) NRMSE = 0.111}\medskip
\end{minipage}
\hfill
\begin{minipage}[b]{.2\linewidth}
  \centering
  \centerline{\includegraphics[height=4cm]{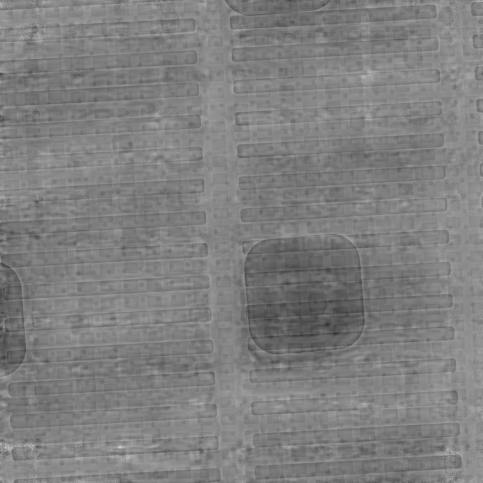}}
  \vspace{0.1cm}
  \centerline{(c) NRMSE = 0.051}\medskip  
\end{minipage}
\hfill
\begin{minipage}[b]{.2\linewidth}
  \centering
  \centerline{\includegraphics[height=4cm]{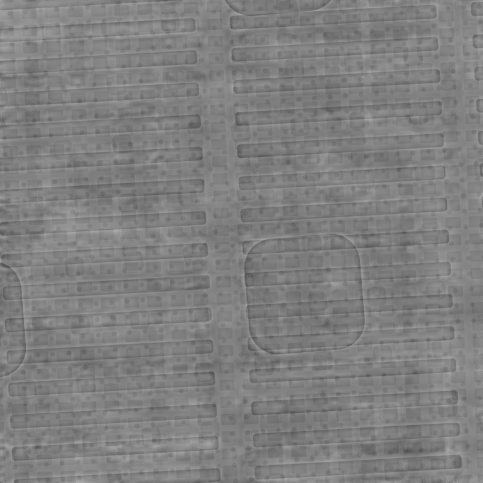}}
  \vspace{0.1cm}
  \centerline{(d) \textbf{NRMSE = 0.037}}\medskip  
\end{minipage}
\hfill
\begin{minipage}[b]{.04\linewidth}
  \centering
  \centerline{\includegraphics[height=4.2cm]{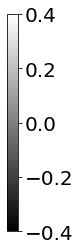}}
  \centerline{ }\medskip  
\end{minipage}
\vspace{0.1cm}
\caption{Top:  {\bf Phase} (in radians) of the reconstructed complex transmittance images in Figure~\ref{fig:ground_truth_imgs} from synthetic data with a small overlap:  $r_{\mathrm{ovlp}} \approx 14.5\%$.  Bottom: Difference between the reconstructed and ground truth phase, with associated NRMSEs in the subcaptions. Note that transmittance phase is typically more important than magnitude in applications. 
}
\label{fig:synthetic_result_phase}
\vspace{-0.2cm}
\end{figure*}

\begin{figure*}[b!]
 \begin{minipage}[b]{.2\linewidth}
  \centering
  \centerline{\includegraphics[height=4cm]{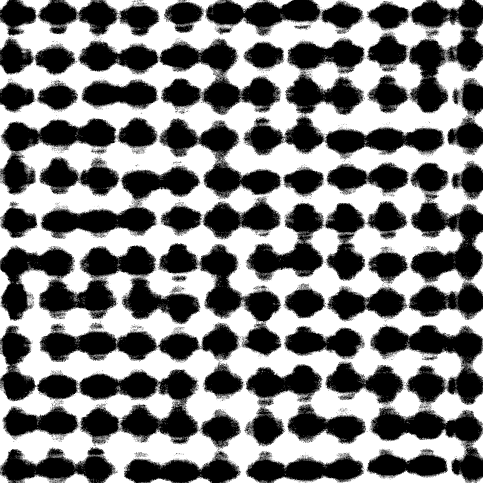}}
  \vspace{0.1cm}
  \centerline{(a) ePIE}\medskip
 \end{minipage}
 \hfill
 \begin{minipage}[b]{.2\linewidth}
  \centering
  \centerline{\includegraphics[height=4cm]{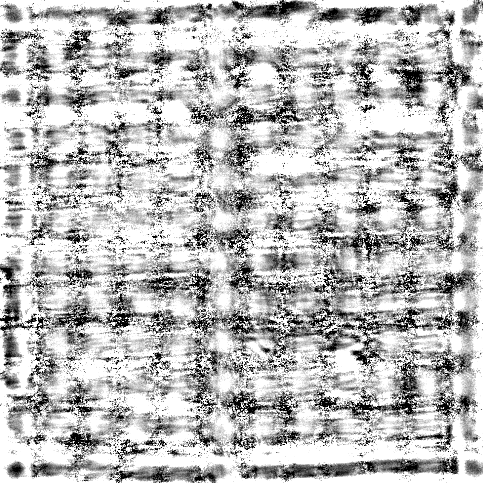}}
  \vspace{0.1cm}
  \centerline{(b) AWF}\medskip
 \end{minipage}
 \hfill
 \begin{minipage}[b]{.2\linewidth}
  \centering
  \centerline{\includegraphics[height=4cm]{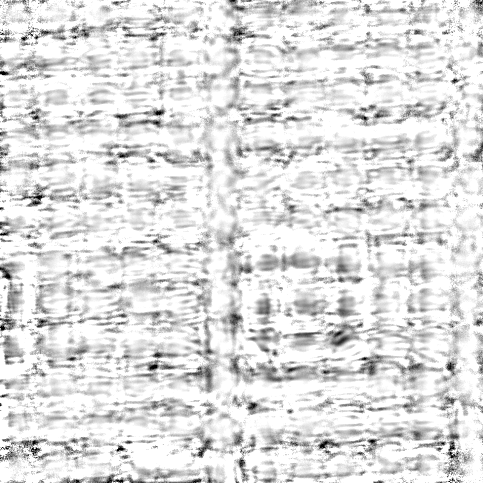}}
  \vspace{0.1cm}
  \centerline{(c) SHARP}\medskip 
 \end{minipage}
 \hfill
 \begin{minipage}[b]{0.2\linewidth}
  \centering
  \centerline{\includegraphics[height=4cm]{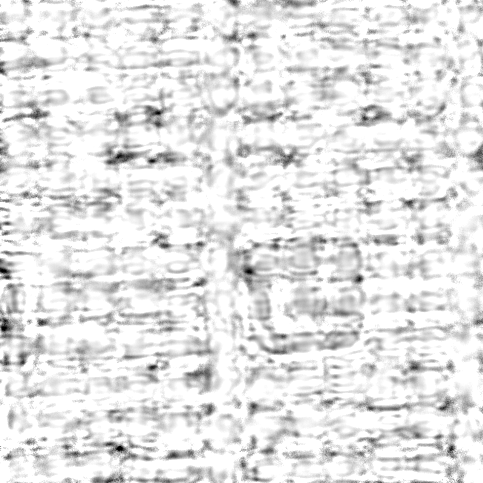}}
  \vspace{0.1cm}
  \centerline{(d) PMACE}\medskip
 \end{minipage}
 \hfill
 \begin{minipage}[b]{.05\linewidth}
  \centering
  \centerline{\includegraphics[height=4.2cm]{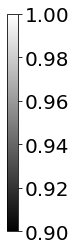}}
  \vspace{0.1cm}
  \centerline{ }\medskip  
 \end{minipage}
 \hfill
 \begin{minipage}[b]{.2\linewidth}
  \centering
  \centerline{\includegraphics[height=4cm]{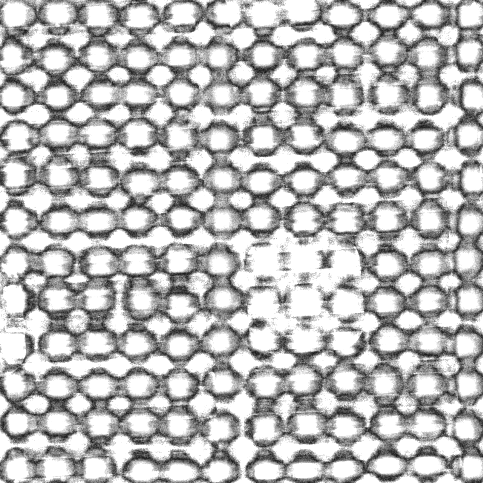}}
  \vspace{0.1cm}
  \centerline{(a) NRMSE = 0.175}\medskip
 \end{minipage}
 \hfill
 \begin{minipage}[b]{.2\linewidth}
  \centering
  \centerline{\includegraphics[height=4cm]{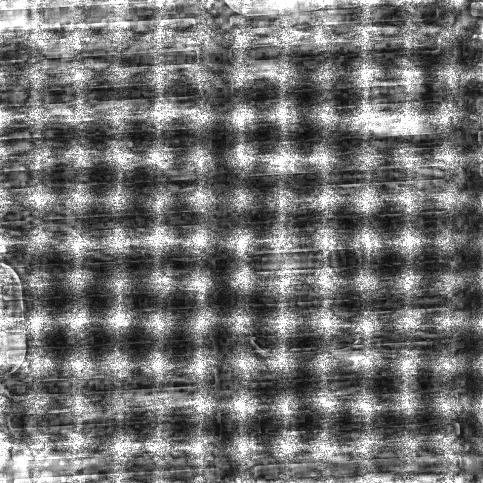}}
  \vspace{0.1cm}
  \centerline{(b) NRMSE = 0.111}\medskip
 \end{minipage}
 \hfill
 \begin{minipage}[b]{.2\linewidth}
  \centering
  \centerline{\includegraphics[height=4cm]{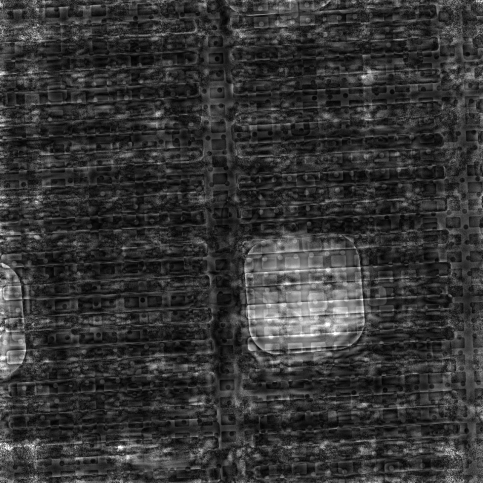}}
  \vspace{0.1cm}
  \centerline{(c) NRMSE = 0.051}\medskip  
 \end{minipage}
 \hfill
 \begin{minipage}[b]{.2\linewidth}
  \centering
  \centerline{\includegraphics[height=4cm]{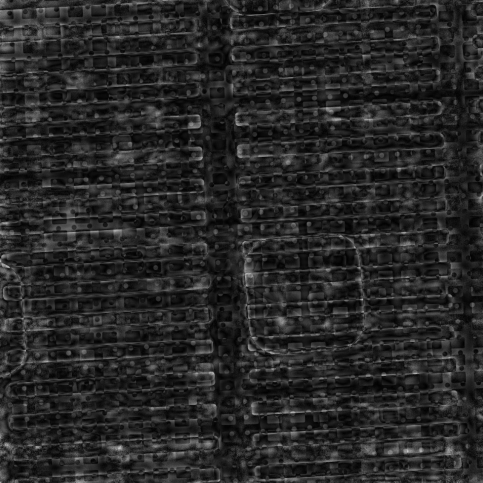}}
  \vspace{0.1cm}
  \centerline{(d) \textbf{ NRMSE = 0.037}}\medskip  
 \end{minipage}
 \hfill
 \begin{minipage}[b]{.04\linewidth}
  \centering
  \centerline{\includegraphics[height=4.2cm]{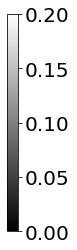}}
  \vspace{0.1cm}
  \centerline{ }\medskip  
 \end{minipage}
 \vspace{0.1cm}
 \caption{Top: {\bf Magnitude} of the reconstructed complex transmittance images in Figure~\ref{fig:ground_truth_imgs} from synthetic data with a small overlap:  $r_{\mathrm{ovlp}} \approx 14.5\%$. Bottom: Amplitudes of error between the complex reconstructions and ground truth, with NRMSEs in the subcaptions.  Note that transmittance phases (Figure~\ref{fig:synthetic_result_phase}) are typically more important in application than magnitudes. }
 \label{fig:synthetic_result_mag}
 \vspace{1cm}
 \end{figure*}

To simulate measurements for a single probe location, we multiply the projection of the complex transmittance image by the complex probe, then take the squared modulus of the Fourier transform to simulate forward propagation. 
We scale to simulate a peak photon rate, add a factor for dark current, then use this as a mean for a sample from a Poisson distribution. 
With $r_{p}$ as the peak photon detection rate and $\lambda$ as the mean dark current, we then have that 
\begin{equation}
\label{eq: poisson_rand_generator}
    \hat{y}_{j} \gets \sqrt{      
\mathrm{Pois}\left ( r_{p} \frac{|\mathcal{F} D x_{j}|^2}{\mathrm{max}_k(\|\mathcal{F} D x_{k}\|_\infty^{2})}  + \lambda \right ) }  \ ,
\end{equation}
where $\|\cdot\|_\infty$ is the maximum magnitude of its argument, $\mathrm{max}_k(\cdot)$ denotes the maximum value over all $k$, and the operations are performed pointwise.
Assuming a photon detector with 14-bit dynamic range and the presence of a half-bit of dark current, $r_{p}=10^4$ and $\lambda = 0.5$ were used to simulate noisy data. 

To simulate data from a full experiment, we first select a nominal distance (in the range of 20 to 76 pixels) between nearest-neighbor probe locations.  We then use a randomly perturbed rectangular grid for the probe sampling pattern with random offsets selected uniformly in $[-5, 5]$ (in units of pixels).
This has the effect of more realistically simulating actual ptychography experiments while also avoiding periodic reconstruction artifacts~\cite{maiden2009improved}.

\subsection{Experimental Methods for Synthetic Data Case}
The key to estimating phase information from intensity measurements is the redundancy inherent in using multiple overlapping illuminations.   
An increase in probe spacing with a fixed probe size will reduce image reconstruction quality~\cite{huang2014optimization}, while using a larger probe will compensate for larger probe spacing~\cite{zhou2020low}. 
Therefore, a measure of illumination overlap is more relevant for reconstruction quality than simple probe spacing.

To quantify this overlap, we define the overlap ratio between adjacent probes indexed by $j$ and $k$ as
\begin{equation}\label{eq:overlap_ratio_adjacent}
    r_{j,k} = \frac{\left \| \; P_{j}^{T} \left | D \right | \odot  P_{k}^{T} \left | D \right |  \; \right \|_{1}}{\left \|  \; \left | D \right | \odot  \left | D \right | \;  \right \|_{1}} \ ,
\end{equation}
where $\odot$ denotes point-wise multiplication and $\left \| \cdot \right \|_{1}$ denotes the $l_1$-norm. We then define the overlap ratio $r_{\mathrm{ovlp}}$ of the full grid scan pattern as the average of the overlap ratios between paired adjacent probes
\begin{equation}\label{eq:overlap_ratio}
   r_{\mathrm{ovlp}} = \frac{1}{N} \sum_{ \left \{ \left (j, k  \right ), \  j\neq k \right \} } r_{j, k} \ ,
\end{equation}
where $N$ denotes the number of pairs of adjacent scan positions.
The overlap ratio is near 1 when adjacent probes almost fully overlap and near 0 when they are almost disjoint. 

To evaluate image reconstruction quality for the synthetic data case, we assume a known complex probe (we will consider joint image-probe estimation in a future paper).  
Given the complex probe profile in Figure~\ref{fig:ground_truth_imgs}, we generated synthetic data using fixed (simulated) probe spacing of 68 pixels with a corresponding overlap ratio of $r_{\mathrm{ovlp}} \approx 14.5\%$. 

All reconstructions are performed using the ePIE~\cite{maiden2009improved}, 
AWF~\cite{xu2018accelerated}, SHARP~\cite{marchesini2016sharp}, and PMACE algorithms.
To provide a fair comparison, we optimized one algorithmic parameter using grid-search for each method: the step size of ePIE, the relaxation parameter of SHARP, and the data-fitting parameter $\alpha$ for PMACE. 
We ran each method for 100 iterations.  
Each iteration for each method involves 2 FFTs per probe location, and these FFTs are the dominant computational cost for each method, so 100 iterations represent essentially the same amount of computation for each method. 

In order to evaluate reconstruction quality, we used the Normalized Root-Mean-Square Error (NRMSE) between each reconstructed image $\hat{x}$ and the ground truth image $x$.
Since ptychography data is invariant to a constant phase shift in the image domain, we computed the NRMSE as  
\begin{equation}
    \begin{aligned}
         NRMSE &= \min_{c} \frac{\| c\hat{x}- x \|}{\| x \|} \ ,
    \end{aligned}
\end{equation}
where $c$ accounts for an unknown phase shift and possible gain.

\subsection{Reconstruction Results on Synthetic Data}

Figures~\ref{fig:synthetic_result_phase} and~\ref{fig:synthetic_result_mag} show the results of the ePIE, AWF, SHARP, and PMACE algorithms.
Figure~\ref{fig:synthetic_result_phase} shows both the phase of the reconstruction in radians and the error in the phase as compared to the ground truth, and it also lists the final NRMSE of the reconstruction.
Figure~\ref{fig:synthetic_result_mag} shows the corresponding images for the magnitude of the reconstruction.

Since the probe spacing is relatively large and the probe overlap is small ($r_{\mathrm{ovlp}} \approx 14.5\%$), ePIE cannot reconstruct the complex transmittance image accurately. 
AWF and SHARP do significantly better, but the difference images indicate that both of these methods introduce artifacts not found in the ground truth.
Overall, PMACE reduces these artifacts and produces the best result both visually and based on the NRMSE.
This is particularly relevant for the phase images, since transmittance phase is typically more important in applications than magnitude.

\begin{figure}[h!]
\centering
\includegraphics[width=3.5in]{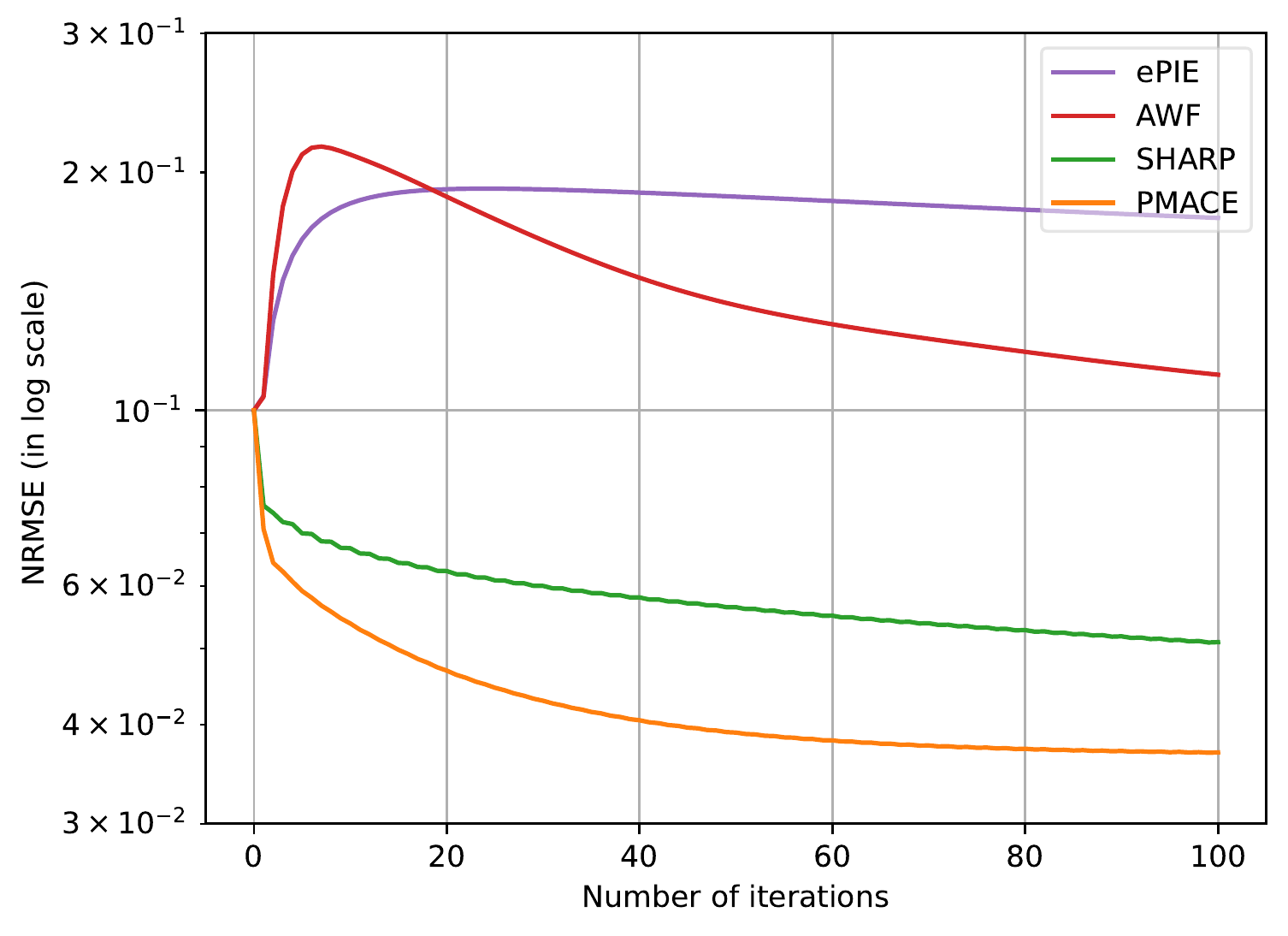}
\caption{NRMSE between reconstruction and ground truth as a function of number of iterations for reconstructions on synthetic noisy data with a small overlap: $ r_{\mathrm{ovlp}} \approx 14.5\% $.}
\label{fig:synthetic_result_convergence_plot}
\end{figure}

Figure~\ref{fig:synthetic_result_convergence_plot} shows the convergence of the NRMSE as a function of the number of iterations.
Consistent with Figures~\ref{fig:synthetic_result_phase} and~\ref{fig:synthetic_result_mag}, the NRMSE associated with ePIE actually increases from the initial condition and does not converge to a reasonable solution. 
We believe this is due to the fact that the problem is very ill-posed due to the small probe overlap.
This is also consistent with other observations on the behavior of ePIE \cite{rodenburg2004phase}.
Similarly, AWF also performs poorly with a large initial error, but then converges slowly to a more reasonable but large final NRMSE. 
SHARP and PMACE have a large initial decrease in error followed by a slower decay, with SHARP leveling off fairly quickly and PMACE continuing to decrease to a good solution.
Ultimately, PMACE achieves the lowest final error.

\subsection{Reconstruction Quality Versus Probe Overlap}

To determine the effect of probe overlap on reconstruction quality, we repeated the experiment above using multiple probe spacing values between 20 and 76 pixels, corresponding to overlap ratios from $73 \%$ to $12 \%$. For each spacing value and each method, we optimized a single algorithmic parameter as described above, and we ran the corresponding method for 100 iterations.  

\begin{figure}[!htbp]
\centering
\includegraphics[width=3.6in]{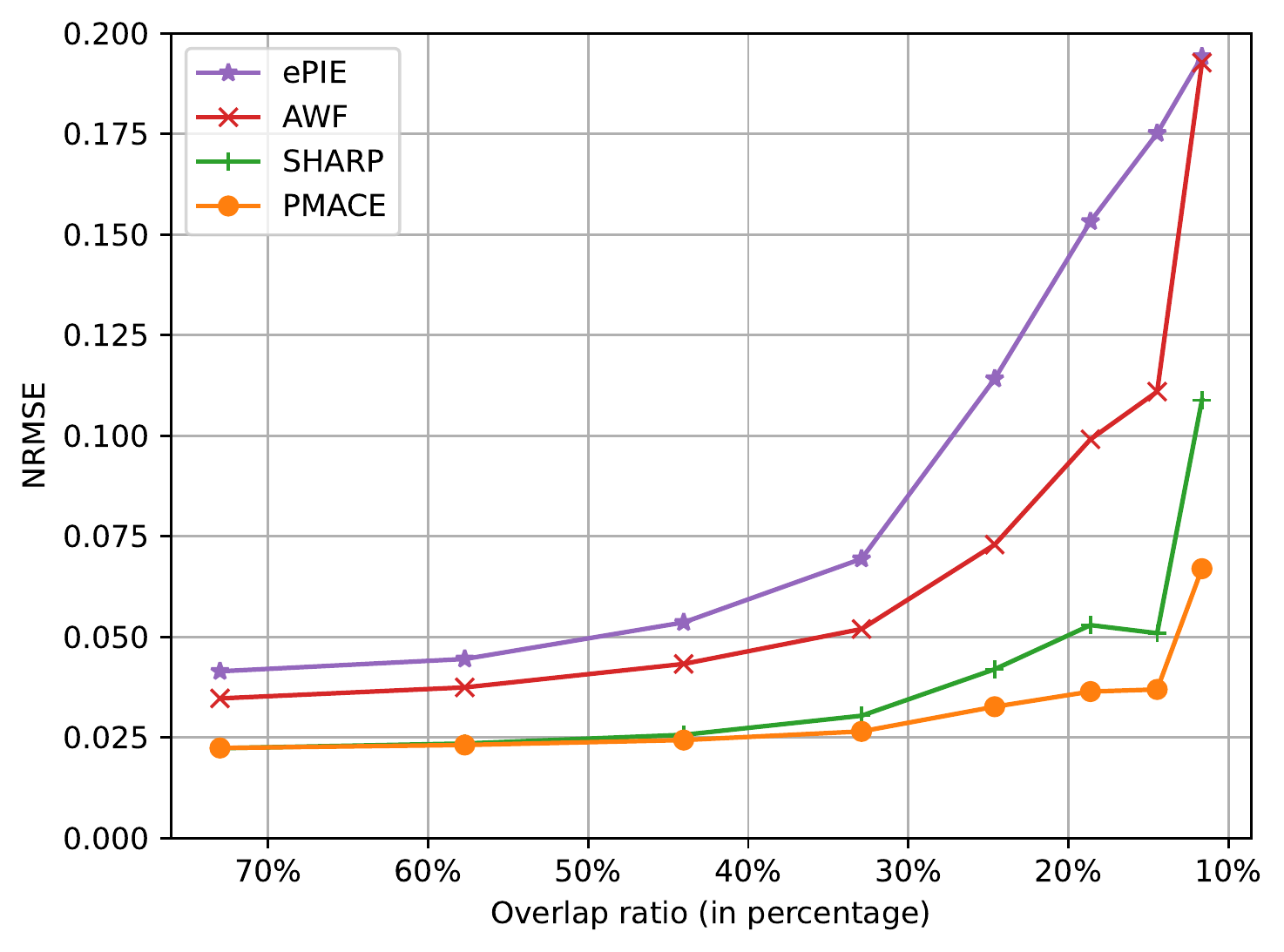}
\caption{NRMSE between reconstruction and ground truth for image of Figure~\ref{fig:ground_truth_imgs} as a function of the overlap ratio for each approach.}
\label{fig: synthetic_result_err_vs_overlap_ratio}
\vspace{-4mm}
\end{figure}

Figure~\ref{fig: synthetic_result_err_vs_overlap_ratio} shows a plot of the final NRMSE of each method versus the overlap ratio. 
For large overlap ratios, all methods yield good image quality, with NRMSE values from about 2.5\% to 4.5\%.  
However, ePIE and AWF underperform relative to SHARP and PMACE.
As the overlap ratio decreases, NRMSE increases, with SHARP and PMACE maintaining good reconstruction quality of about 5\% or less as the overlap ratio decreases to about 15\%.
PMACE achieves an NRMSE below 7\% even with an overlap ratio $r_{\mathrm{ovlp}}$ below 12\%, while the SHARP reconstruction has an NRMSE above 10\% at this overlap ratio.  ePIE and AWF are not competitive at small overlap ratios.

\begin{figure}[!htbp]
\centering
\includegraphics[width=3.8in]{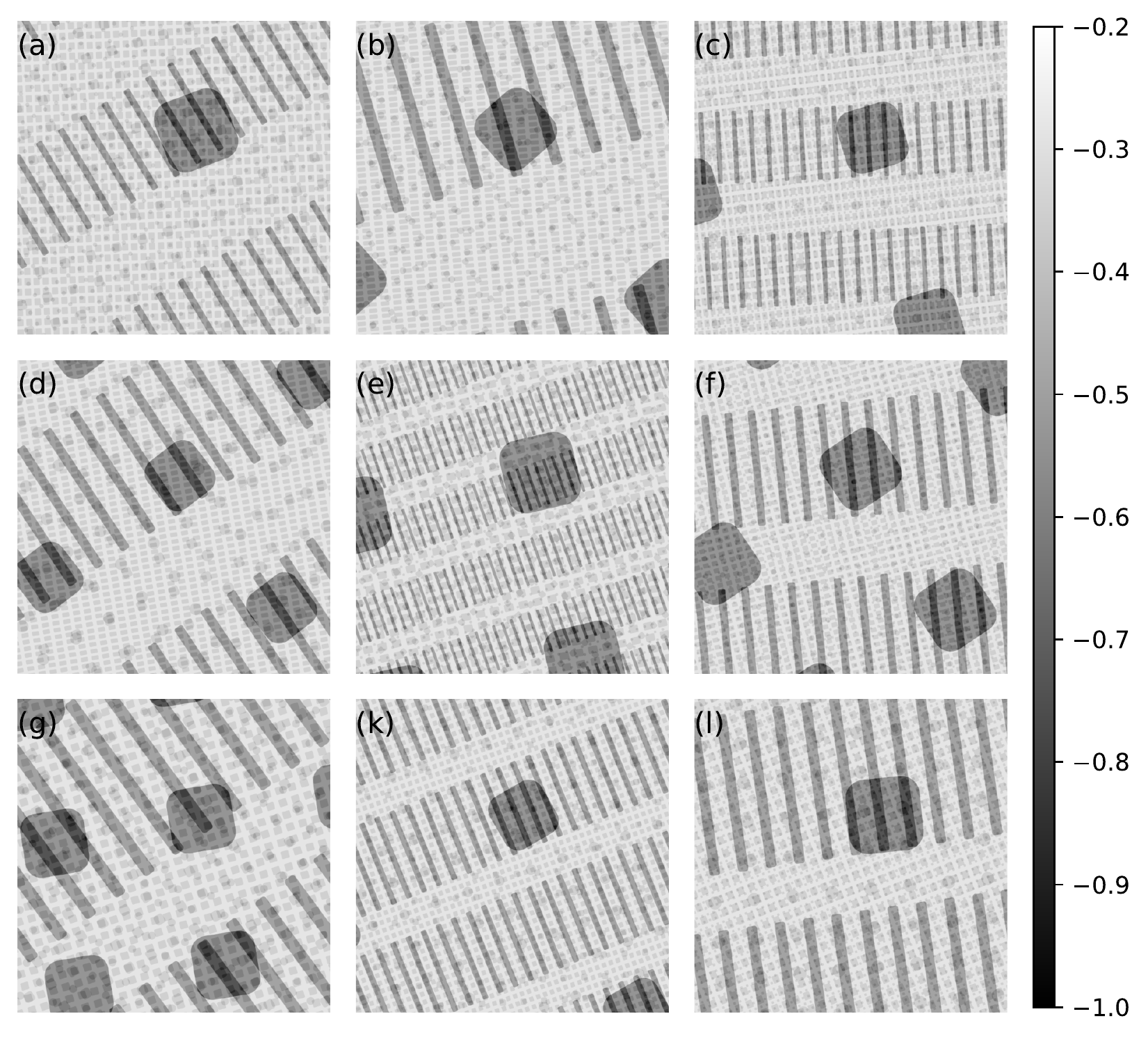}
\caption{Phase of the 9 additional complex transmittance images used as ground truth in simulation experiments of Figure~~\ref{fig: avg_nrmse_vs_overlap_ratio}.}
\label{fig: sample_imgs}
\end{figure}

\begin{figure}[t]
\centering
\includegraphics[width=3.6in]{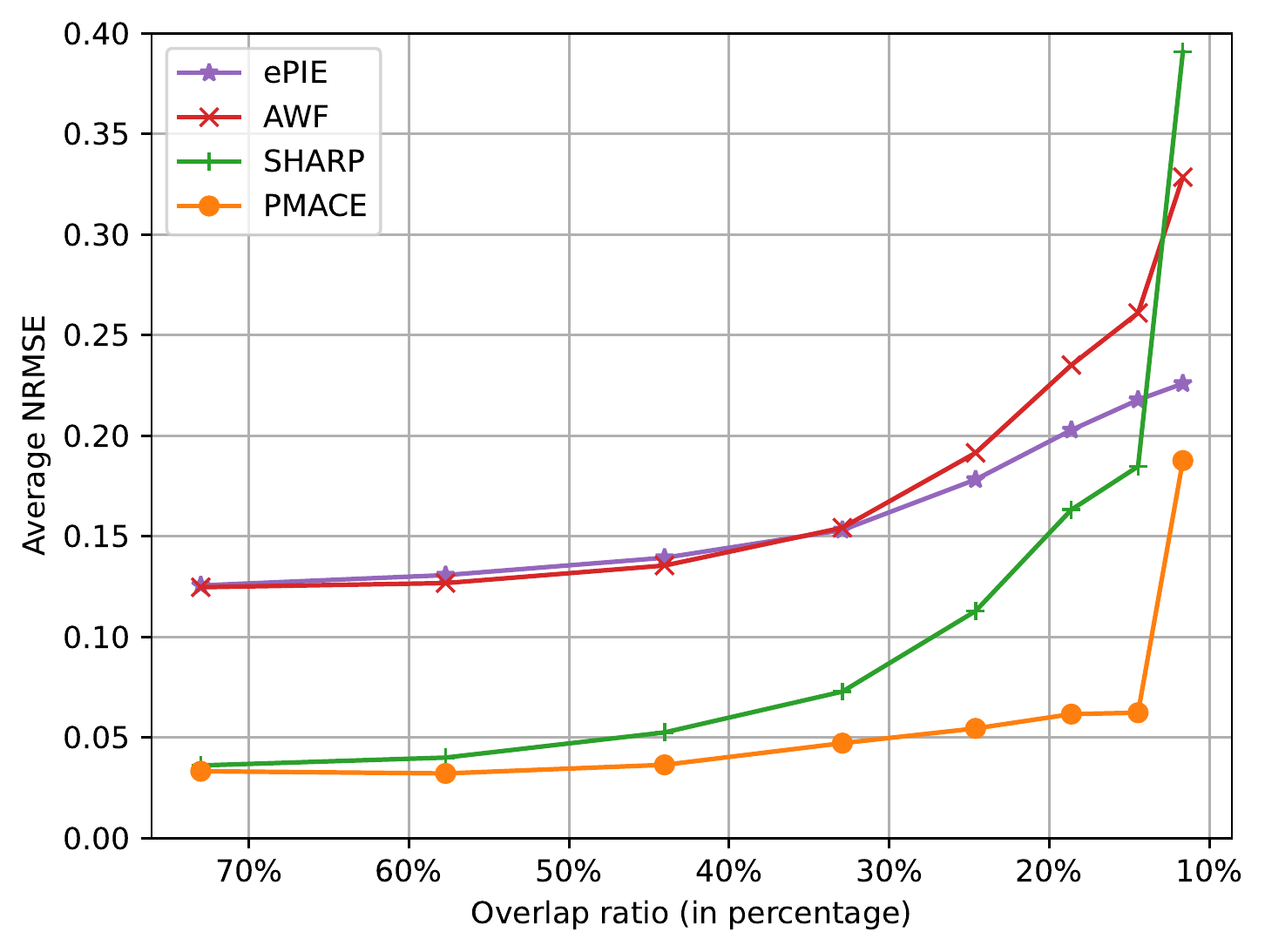}
\caption{Average NRMSE between reconstructions and ground truth as a function of the overlap ratio for each approach.}
\label{fig: avg_nrmse_vs_overlap_ratio}
\end{figure}

Figure~\ref{fig: avg_nrmse_vs_overlap_ratio} shows a comparison of the algorithms for a set of 9 distinct $800\times 800$ ground truth images shown in Figure~\ref{fig: sample_imgs}.
Again, this result demonstrates that PMACE maintains the best reconstruction quality when the overlap percentage is low, and that SHARP and PMACE are significantly better than ePIE and AWF.

\subsection{Reconstruction Results on Measured Data}
In this subsection, we evaluate PMACE on measured data and report reconstruction results.  
We used the gold balls data set~\cite{marchesini2017ptychography} consisting of 800 ptychographic measurements, which was generated using X-rays with a photon energy of 1 keV at beamline 5.3.2 of the Advanced Light Source. The scans were performed in a $20 \times 40$ grid with scan spacing of 30 nm and recorded by a detector placed 0.112 m downstream from the test sample of gold balls.
This data set also includes 20 separately acquired dark images in addition to measured data.

For this experiment, the overlap ratio is approximately $84\%$, which is higher than any value used in our simulated data. This large overlap ratio and the results of Figure~\ref{fig: synthetic_result_err_vs_overlap_ratio} indicate that any of the methods should be able to produce a reasonably accurate reconstruction.  

We preprocessed the raw data by subtracting from each measurement the average of 20 dark scans.  We also removed 6 abnormal measurements with high deviations in the data. Then we centered and cropped each measurement from $621 \times 621$ pixels to $512 \times 512$ pixels. 
To suppress noise, each of these diffraction measurements was multiplied by a 2D Tukey window that was generated by rotating a 1D Tukey window with shape parameter of 0.5. The resulting data set contains preprocessed 794 diffraction measurements. 

\begin{figure}[!htbp]
\vspace{-0.4cm}
\begin{minipage}[htb]{0.45\linewidth}
  \centering
   \centerline{\includegraphics[height=4cm]{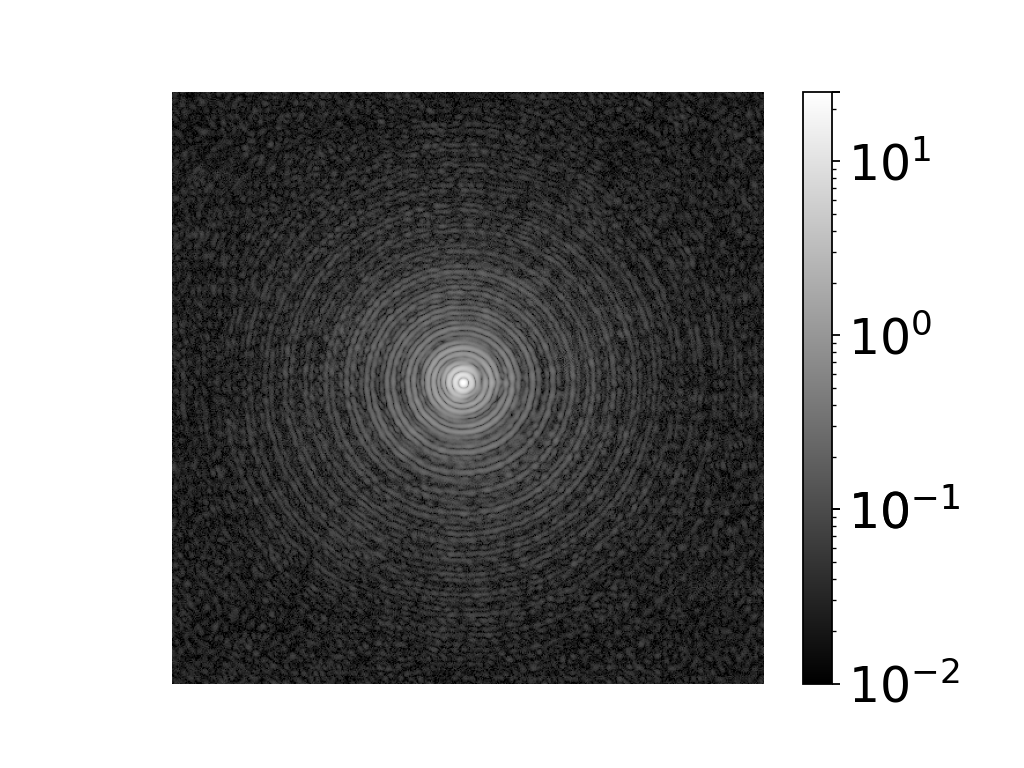}}
  \vspace{-0.4cm}
  \centerline{(a)} \smallskip
\end{minipage}
\hspace{0.1cm}
\begin{minipage}[htb]{0.45\linewidth}
  \centering
  \centerline{\includegraphics[height=4cm]{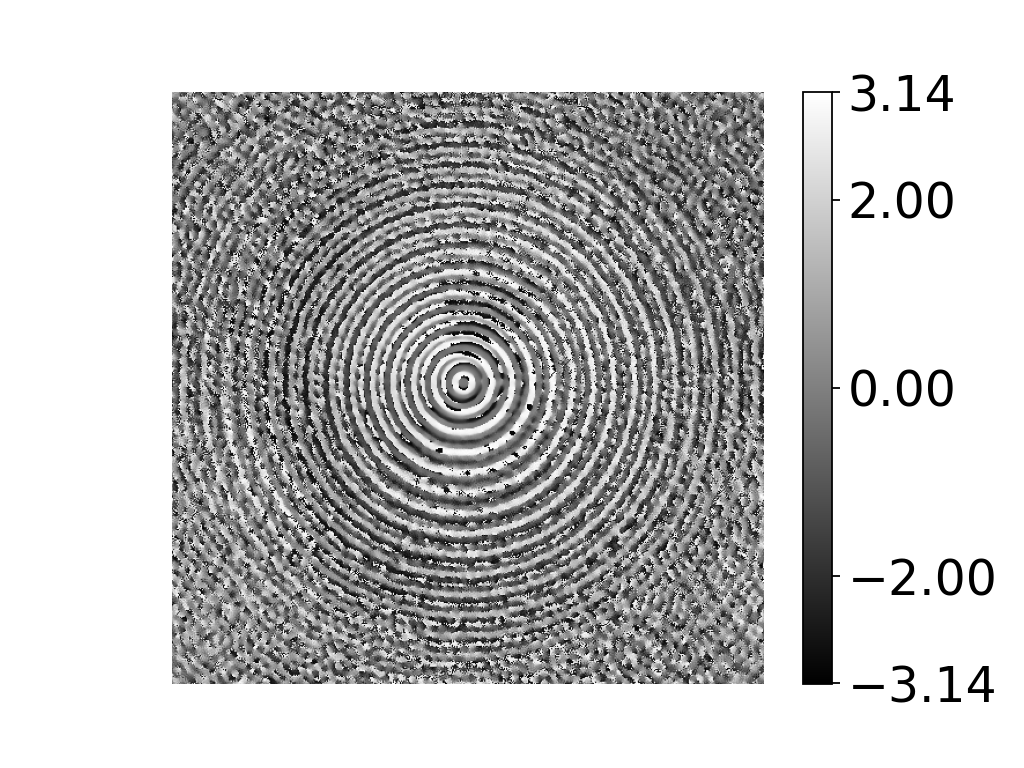}}
  \vspace{-0.4cm}
  \centerline{(b)} \smallskip
\end{minipage}
\hfill
\caption{
The (a) magnitude and (b) phase of the probe used in all experiments on the measured data. The probe was estimated from the data using a variation of the ePIE algorithm.
}
\label{fig:ref_probe}
\end{figure}

Since the goal of this research is to compare the quality and speed of different ptychographic reconstruction algorithms, we used a fixed probe function for all algorithms in our comparisons.
This fixed probe, shown in Figure~\ref{fig:ref_probe}, was estimated from experimental data using a variation of the ePIE algorithm.
We note that probe estimation is important in practice and has been incorporated into other algorithms, but we leave the problem of PMACE probe estimation as a topic for future research.

As in the case of simulated data, we varied one algorithmic parameter for each method and selected the parameter that produced the highest image contrast and fewest image artifacts. We ran each method for 100 iterations. 
As noted above, the fit to data is unchanged by a constant phase shift in the reconstruction.  To provide a uniform visual comparison, we shift the phase of each reconstruction to achieve zero-mean phase in the central rectangle of $105 \times 65$ pixels, which was chosen since the centers of the probe illuminations are located in this region.  

\begin{figure*}[!htbp]
\begin{minipage}[b]{.2\linewidth}
  \centering
  \centerline{\includegraphics[height=4cm]{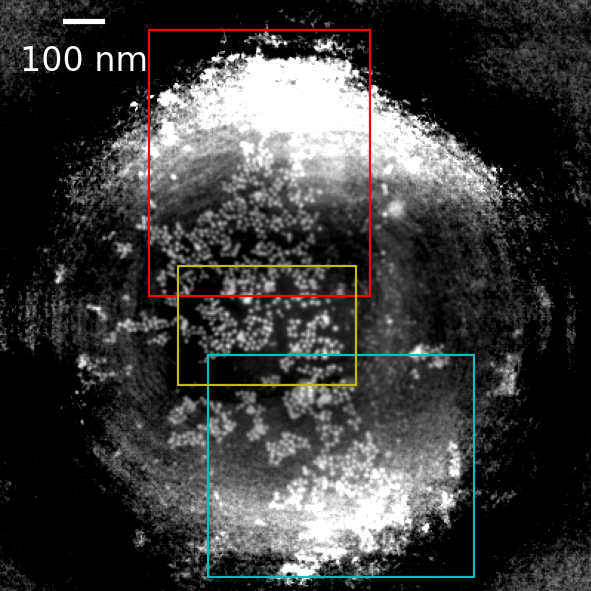}}
  \vspace{0.1cm}
  \centerline{(a) ePIE}\medskip
\end{minipage}
\hfill
\begin{minipage}[b]{.2\linewidth}
  \centering
  \centerline{\includegraphics[height=4cm]{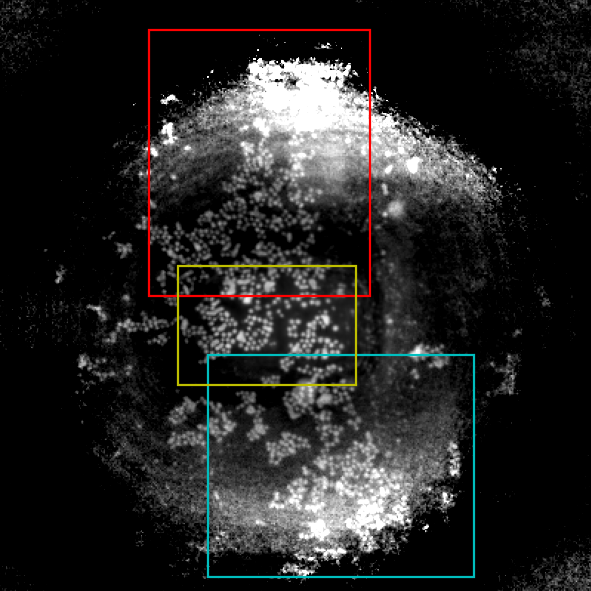}}
  \vspace{0.1cm}
  \centerline{(b) AWF}\medskip
\end{minipage}
\hfill
\begin{minipage}[b]{.2\linewidth}
  \centering
  \centerline{\includegraphics[height=4cm]{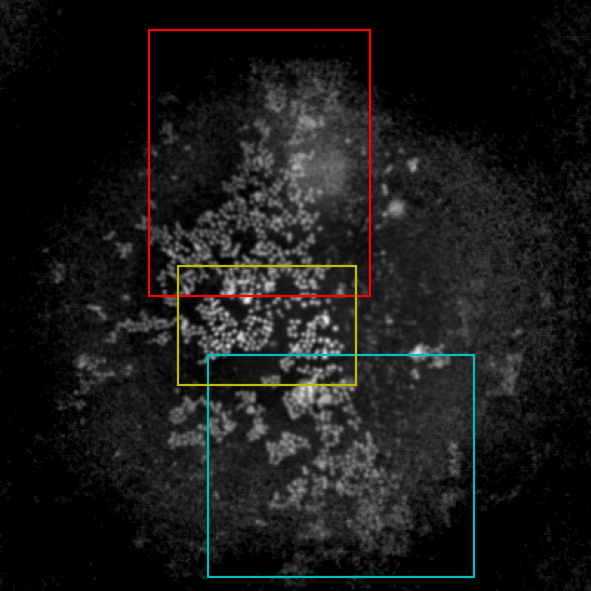}}
  \vspace{0.1cm}
  \centerline{(c) SHARP}\medskip  
\end{minipage}
\hfill
\begin{minipage}[b]{.2\linewidth}
  \centering
  \centerline{\includegraphics[height=4cm]{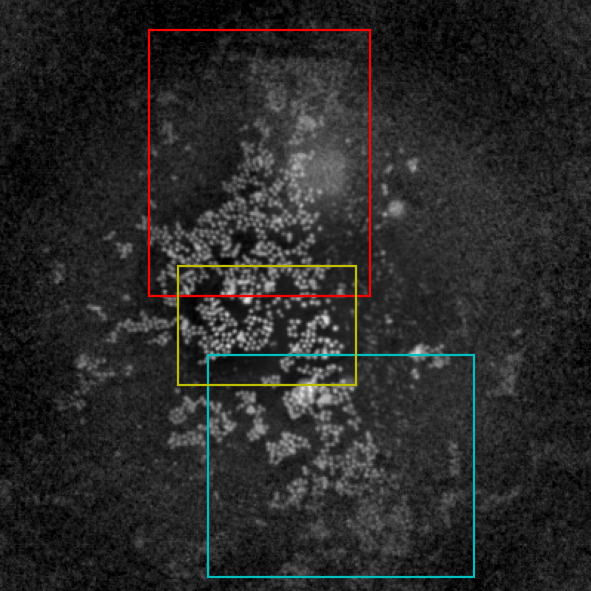}}
  \vspace{0.1cm}
  \centerline{(d) PMACE}\medskip  
\end{minipage}
\hfill
\begin{minipage}[b]{.04\linewidth}
  \centering
  \centerline{\includegraphics[height=4.2cm]{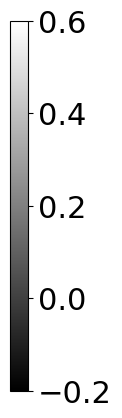}}
  \centerline{}\medskip
\end{minipage}
\hfill
\begin{minipage}[b]{.2\linewidth}
  \centering
  \vspace{-0.4cm}
  \centerline{\includegraphics[height=4.2cm]{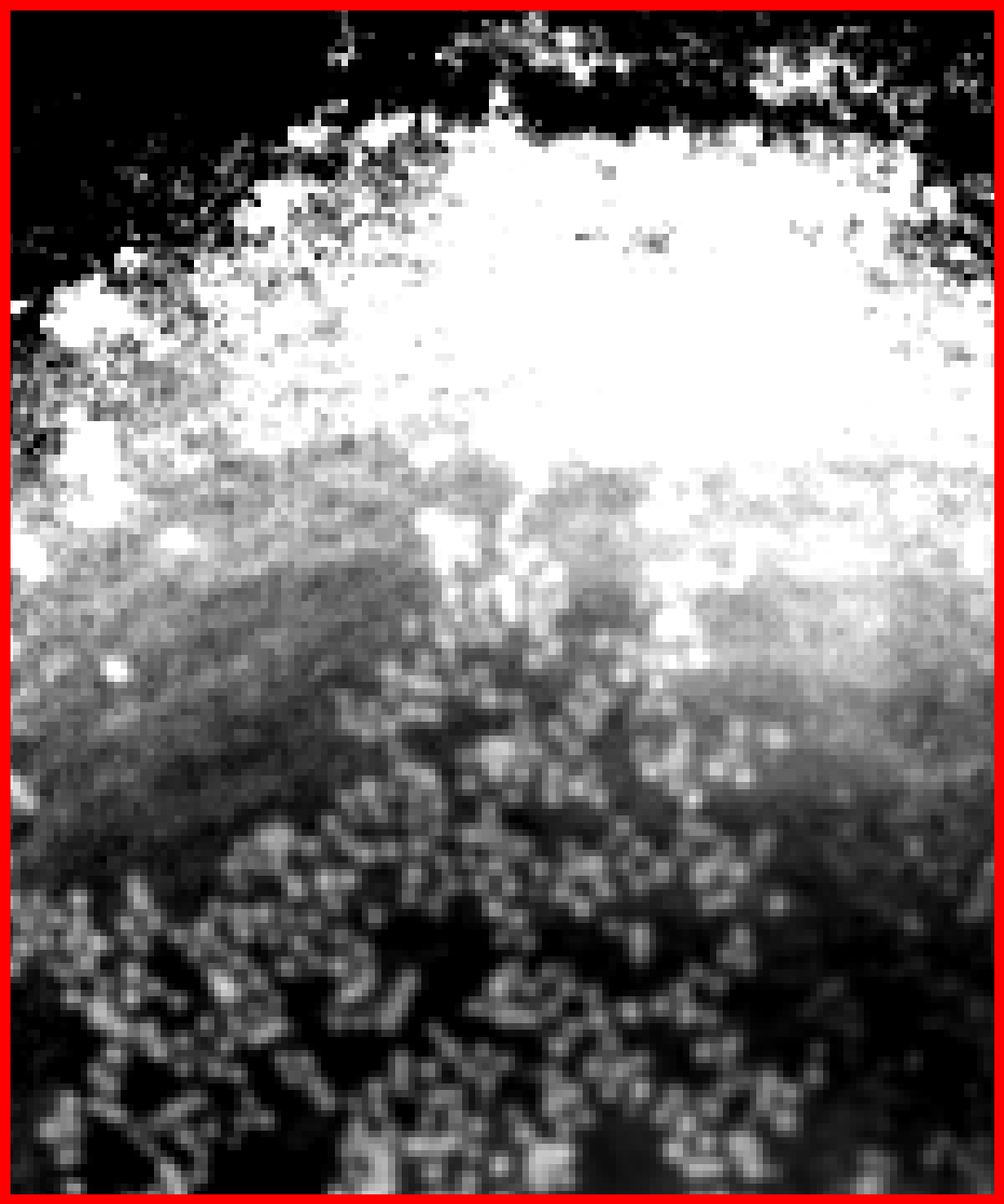}}
  \vspace{0.1cm}
\end{minipage}
\hfill
\begin{minipage}[b]{.2\linewidth}
  \centering
  \vspace{-0.4cm}
  \centerline{\includegraphics[height=4.2cm]{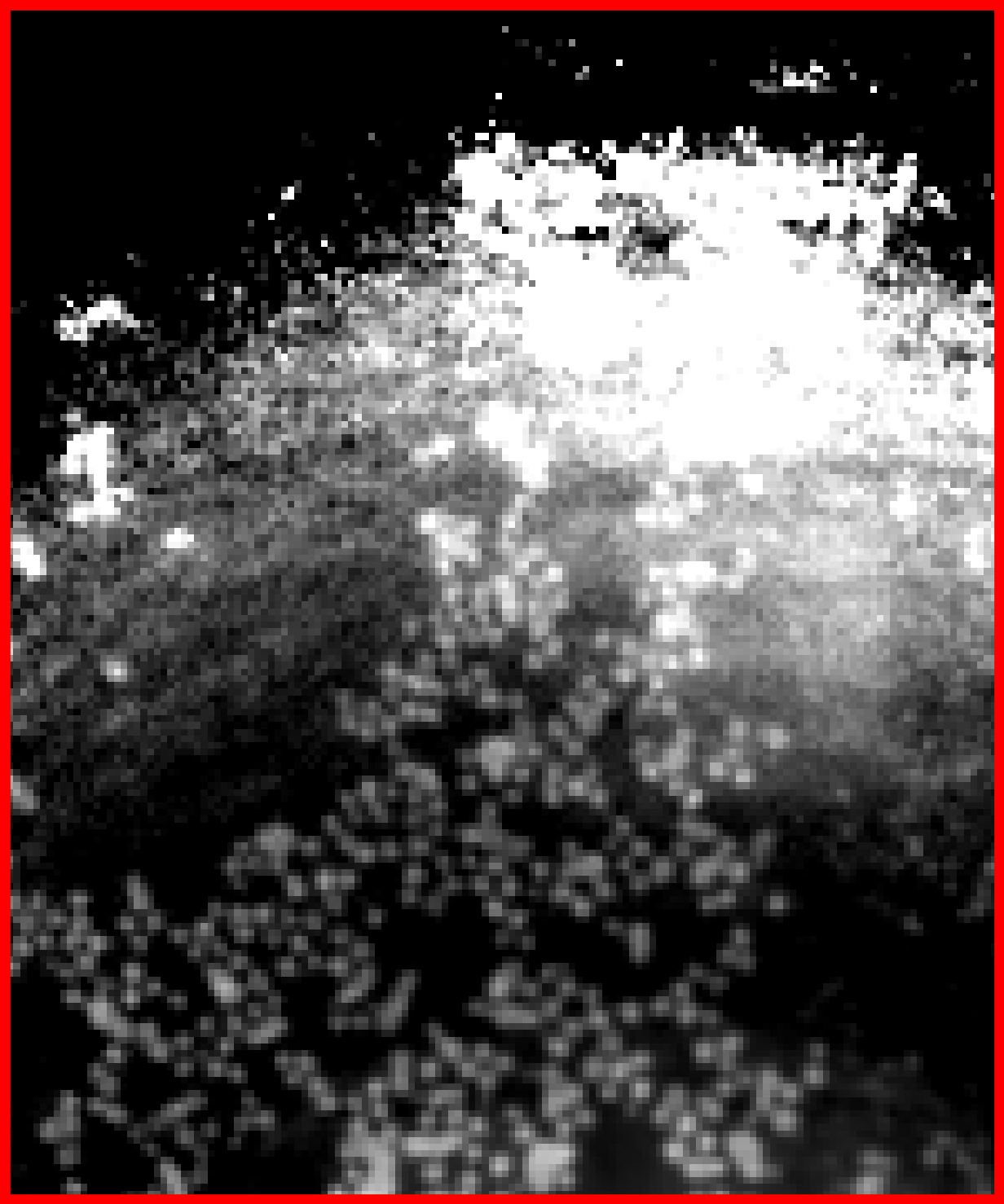}}
  \vspace{0.1cm}
\end{minipage}
\hfill
\begin{minipage}[b]{.2\linewidth}
  \centering
  \vspace{-0.4cm}
  \centerline{\includegraphics[height=4.2cm]{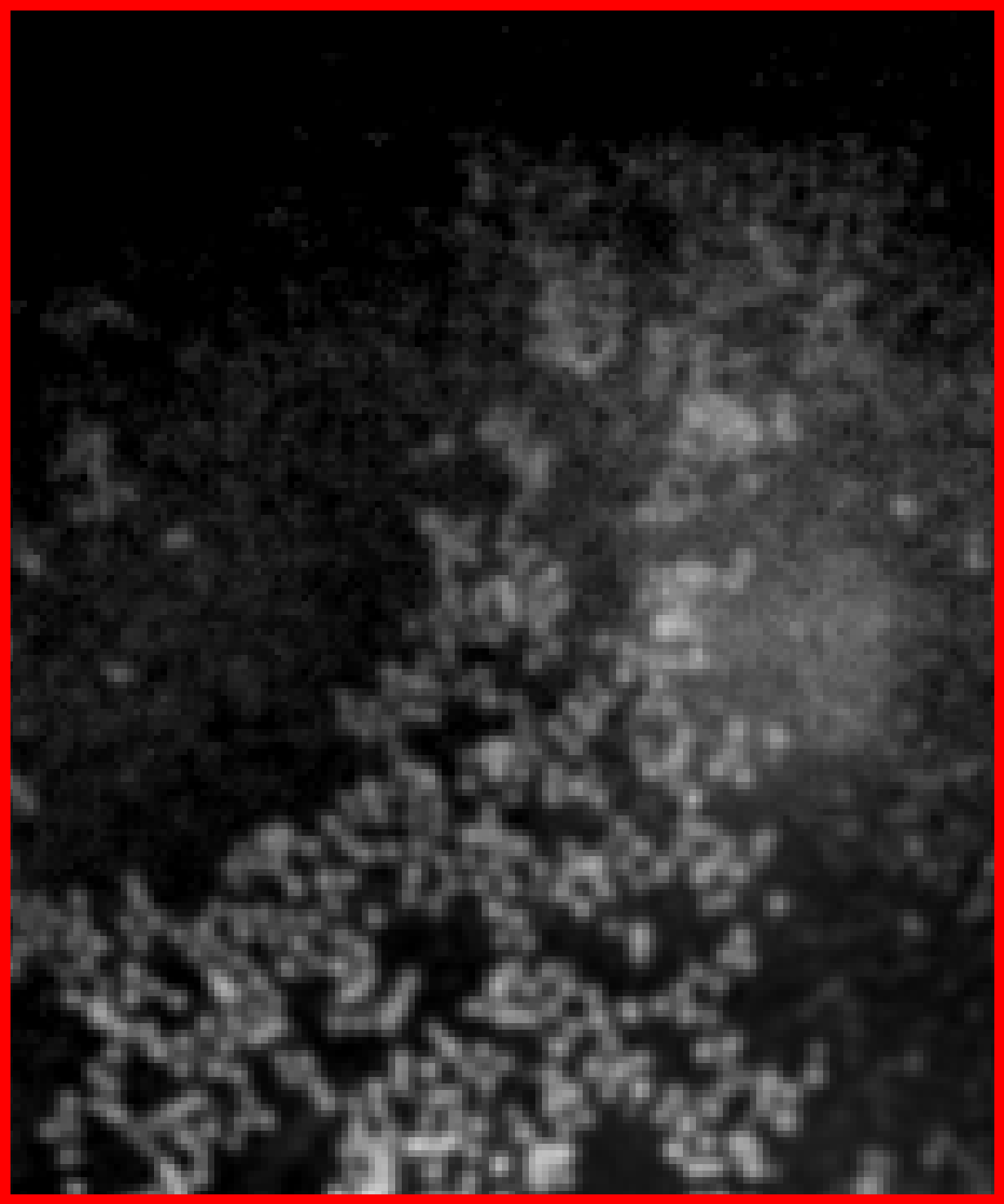}}
  \vspace{0.1cm}
\end{minipage}
\hfill
\begin{minipage}[b]{.2\linewidth}
  \centering
  \vspace{-0.4cm}
  \centerline{\includegraphics[height=4.2cm]{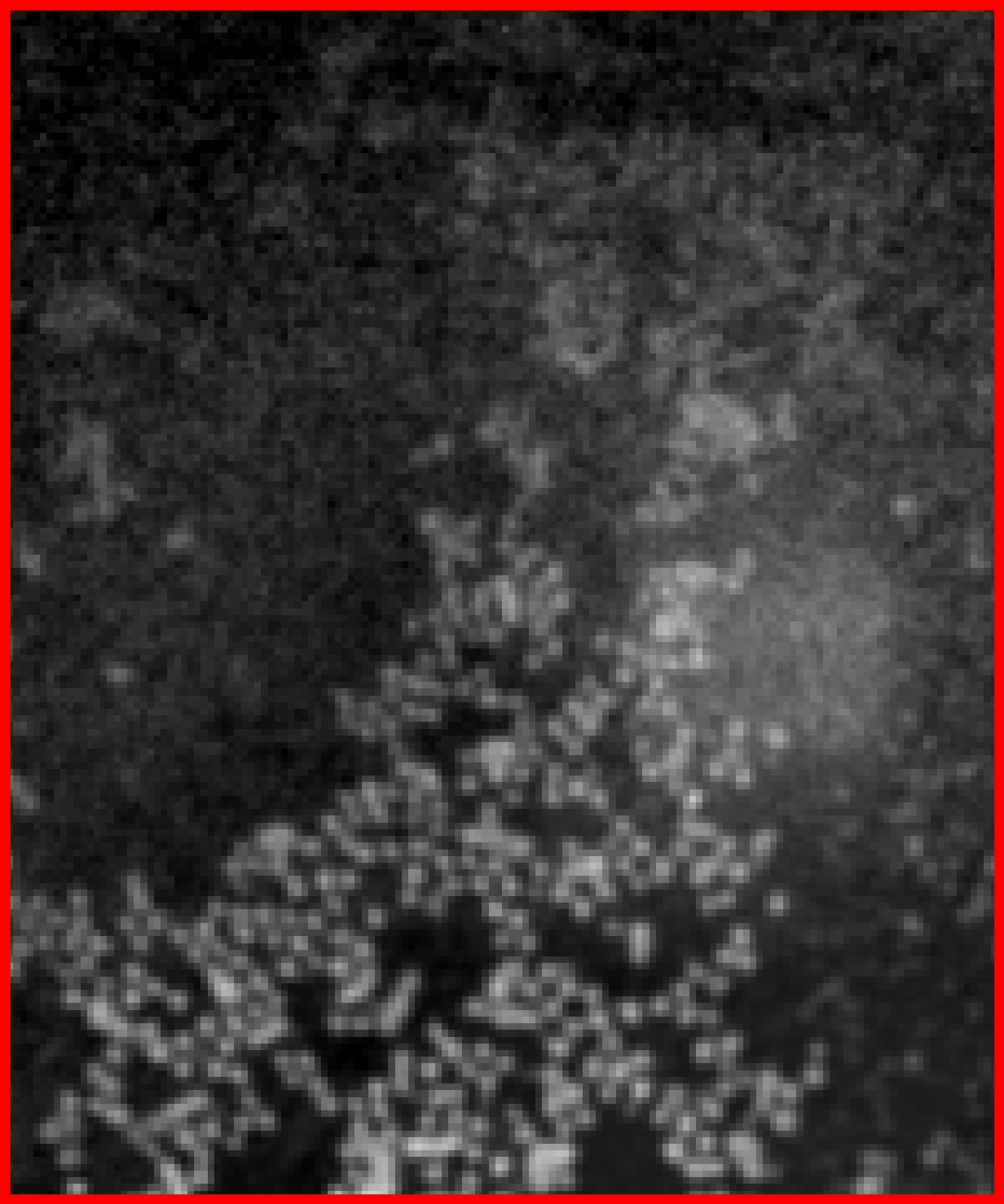}}
  \vspace{0.1cm}
\end{minipage}
\hfill
\begin{minipage}[b]{.04\linewidth}
  \centering
  \vspace{-0.4cm}
  \centerline{\includegraphics[height=4.2cm]{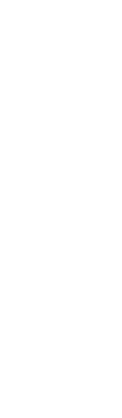}}
  \centerline{}\medskip
\end{minipage}
\hfill
\begin{minipage}[b]{.2\linewidth}
  \centering
  \centerline{\includegraphics[height=2.2cm]{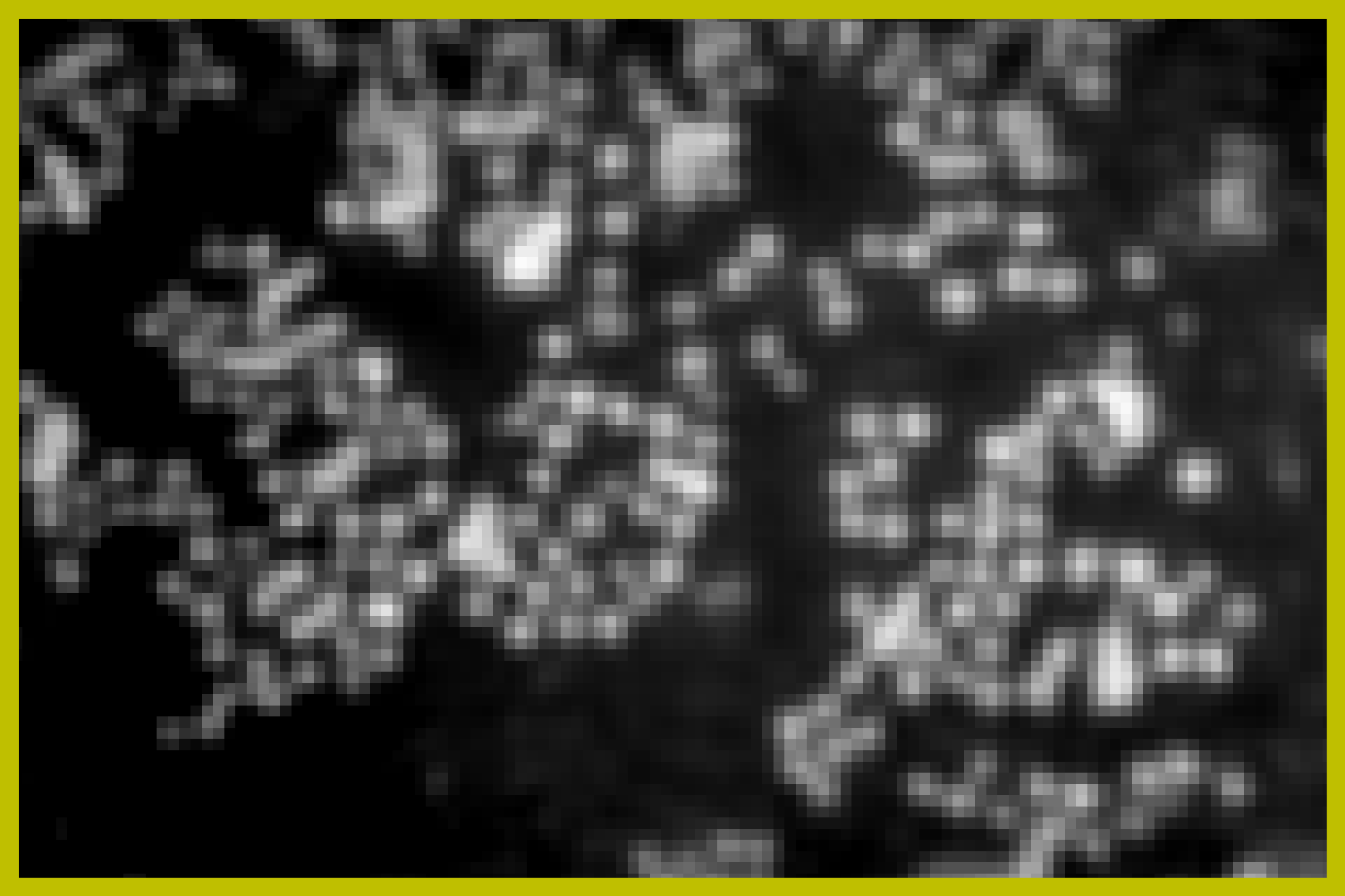}}
  \vspace{0.36cm}
\end{minipage}
\hfill
\begin{minipage}[b]{.2\linewidth}
  \centering
  \centerline{\includegraphics[height=2.2cm]{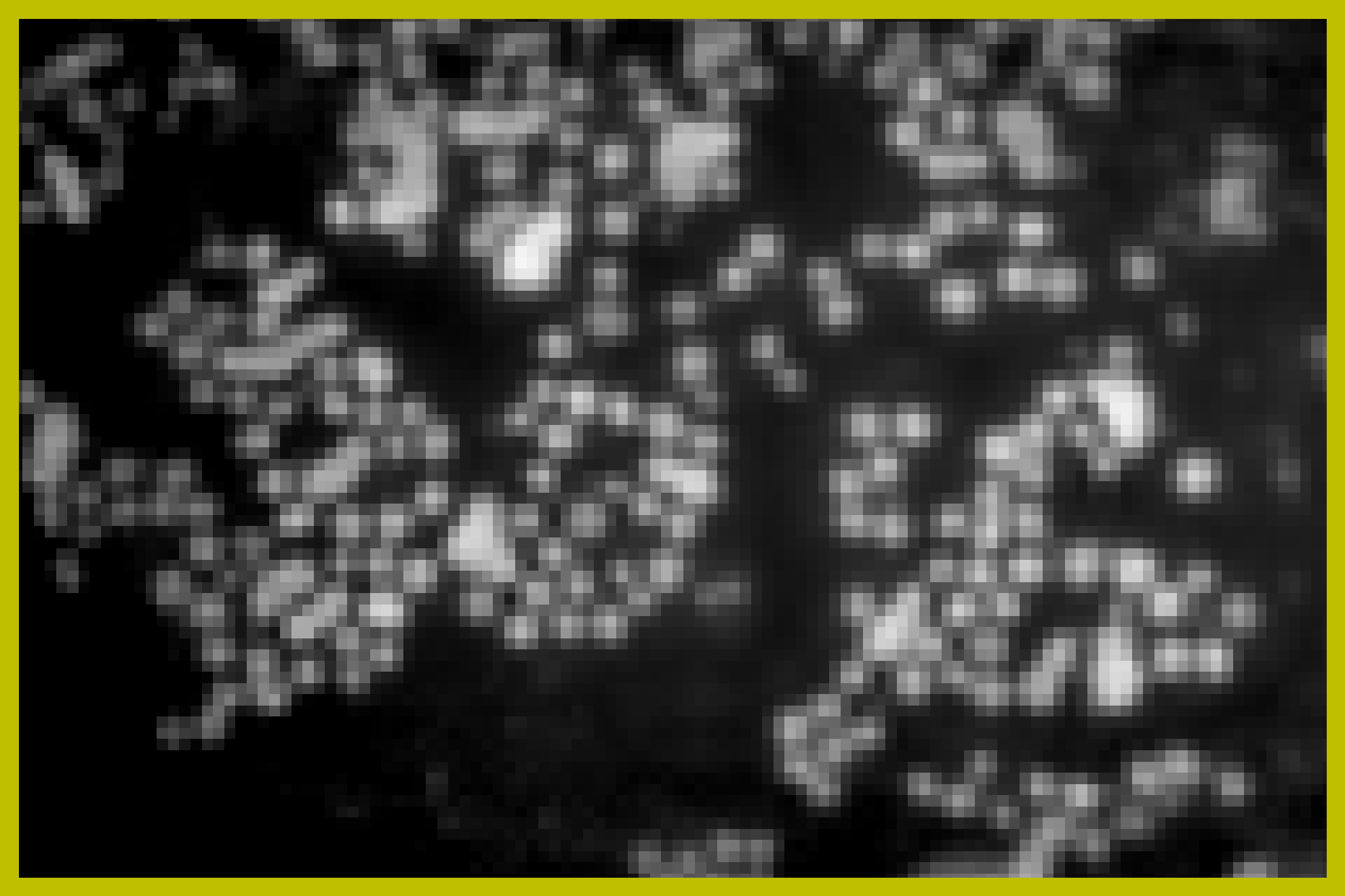}}
  \vspace{0.36cm}
\end{minipage}
\hfill
\begin{minipage}[b]{.2\linewidth}
  \centering
  \centerline{\includegraphics[height=2.2cm]{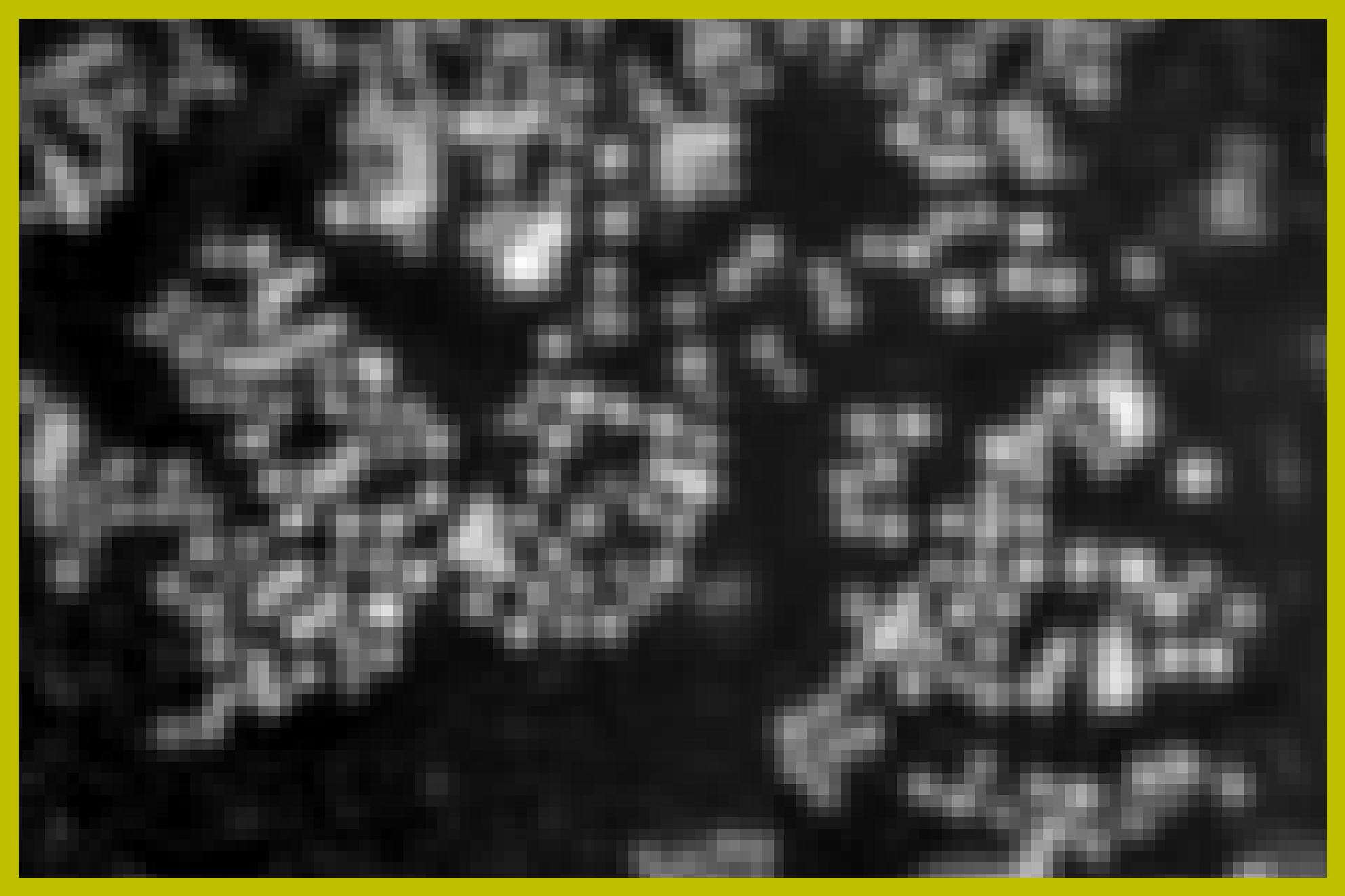}}
  \vspace{0.36cm}
\end{minipage}
\hfill
\begin{minipage}[b]{.2\linewidth}
  \centering
  \centerline{\includegraphics[height=2.2cm]{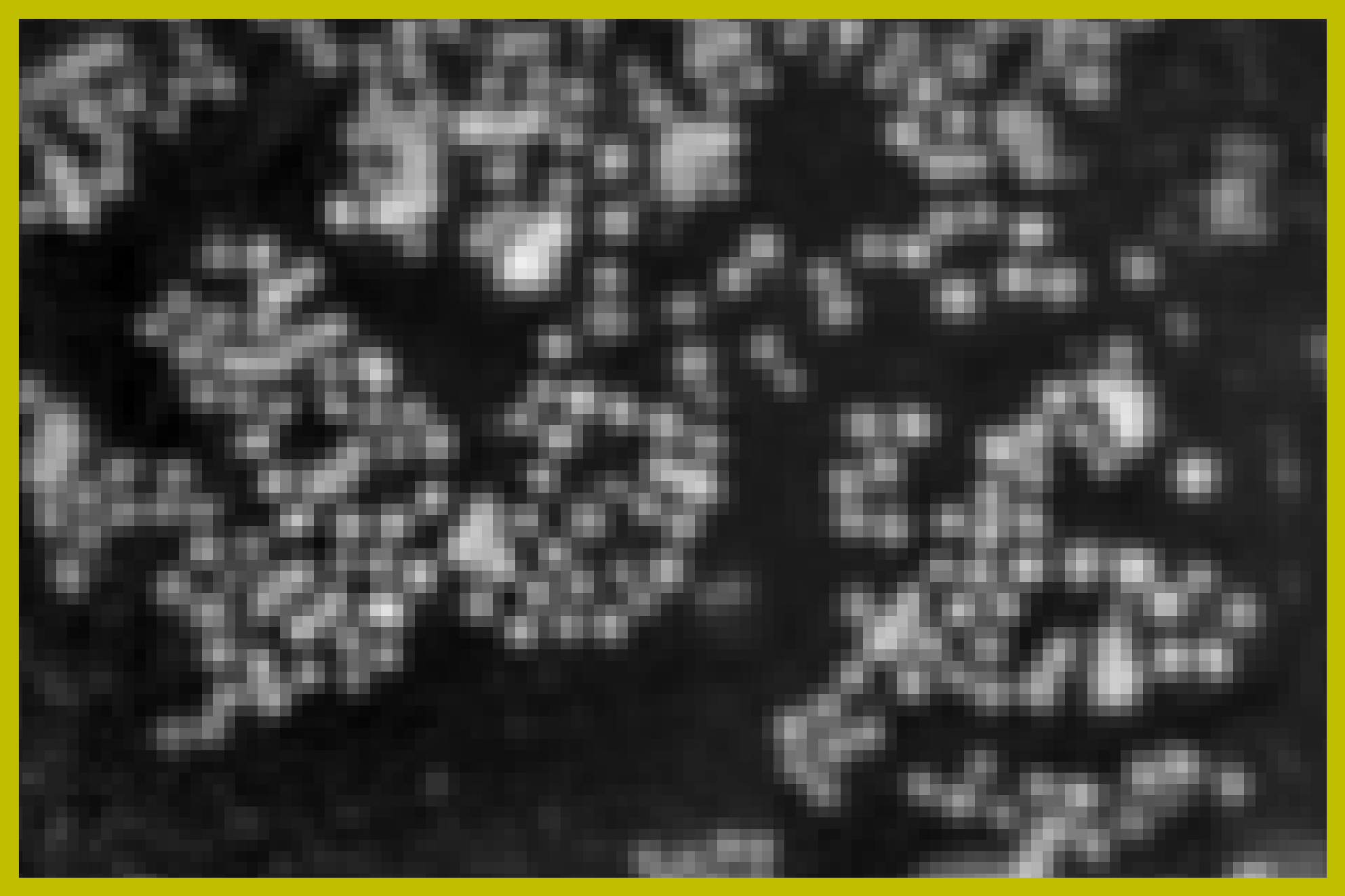}}
  \vspace{0.36cm}
\end{minipage}
\hfill
\begin{minipage}[b]{.04\linewidth}
  \centering
  \centerline{\includegraphics[height=2.2cm]{figs/GoldBalls/colorbar_phase_mask.jpg}}
  \centerline{}\medskip
\end{minipage}
\hfill
\begin{minipage}[b]{.2\linewidth}
  \centering
  \centerline{\includegraphics[height=3.4cm]{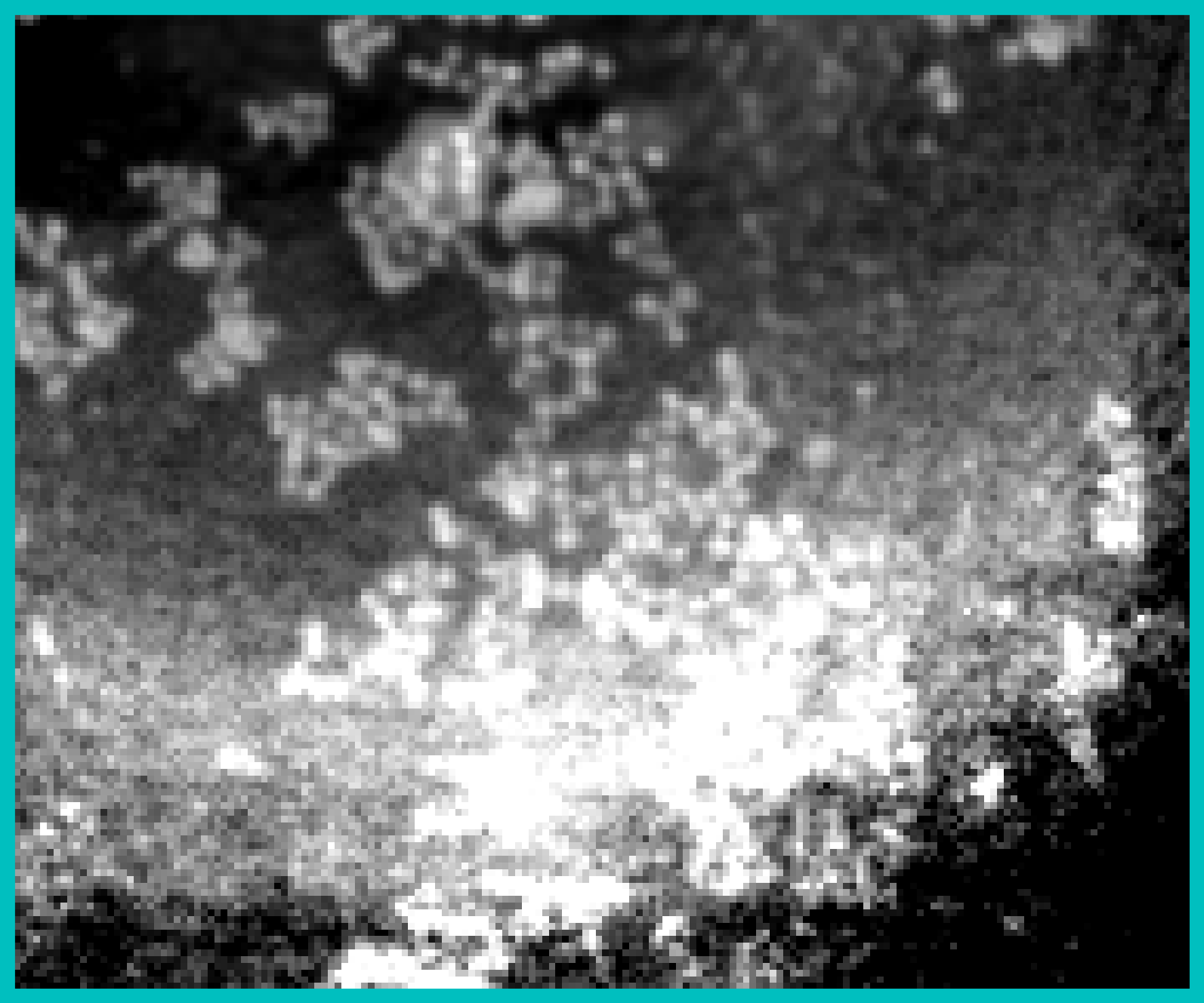}}
  \vspace{0.1cm}
  \centerline{(a) ePIE}\medskip
\end{minipage}
\hfill
\begin{minipage}[b]{.2\linewidth}
  \centering
  \centerline{\includegraphics[height=3.4cm]{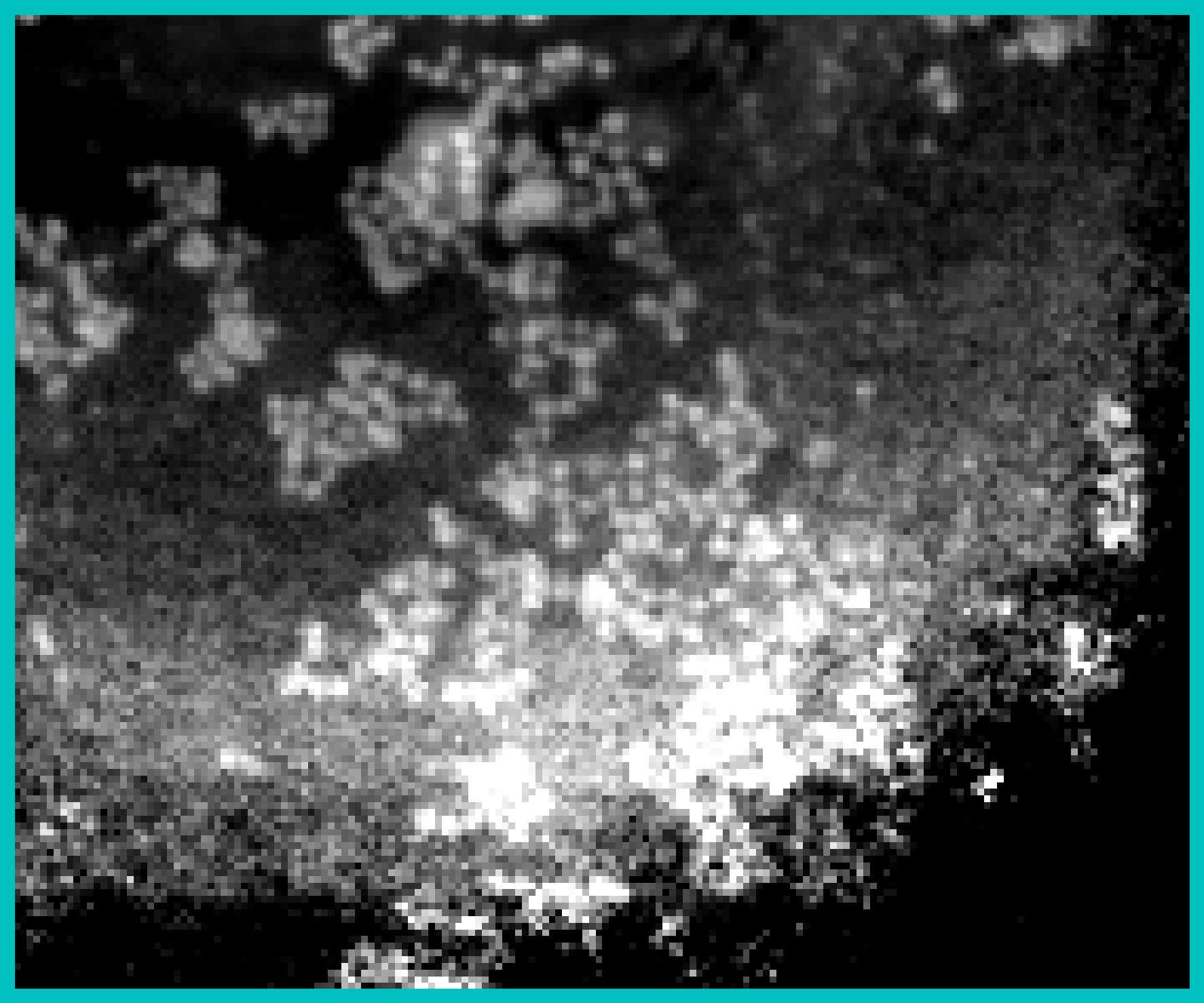}}
  \vspace{0.1cm}
  \centerline{(b) AWF}\medskip
\end{minipage}
\hfill
\begin{minipage}[b]{.2\linewidth}
  \centering
  \centerline{\includegraphics[height=3.4cm]{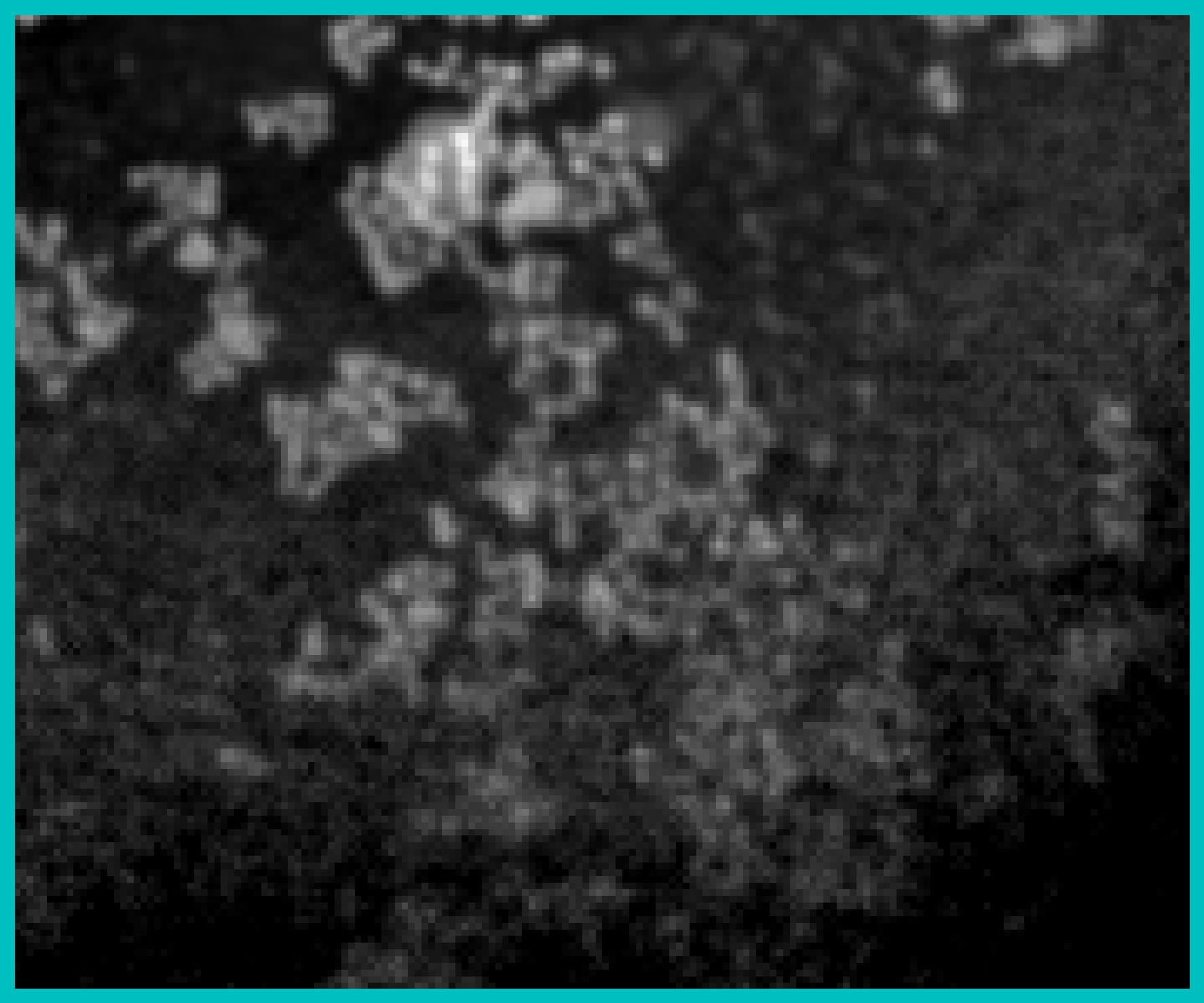}}
  \vspace{0.1cm}
  \centerline{(c) SHARP}\medskip  
\end{minipage}
\hfill
\begin{minipage}[b]{.2\linewidth}
  \centering
  \centerline{\includegraphics[height=3.4cm]{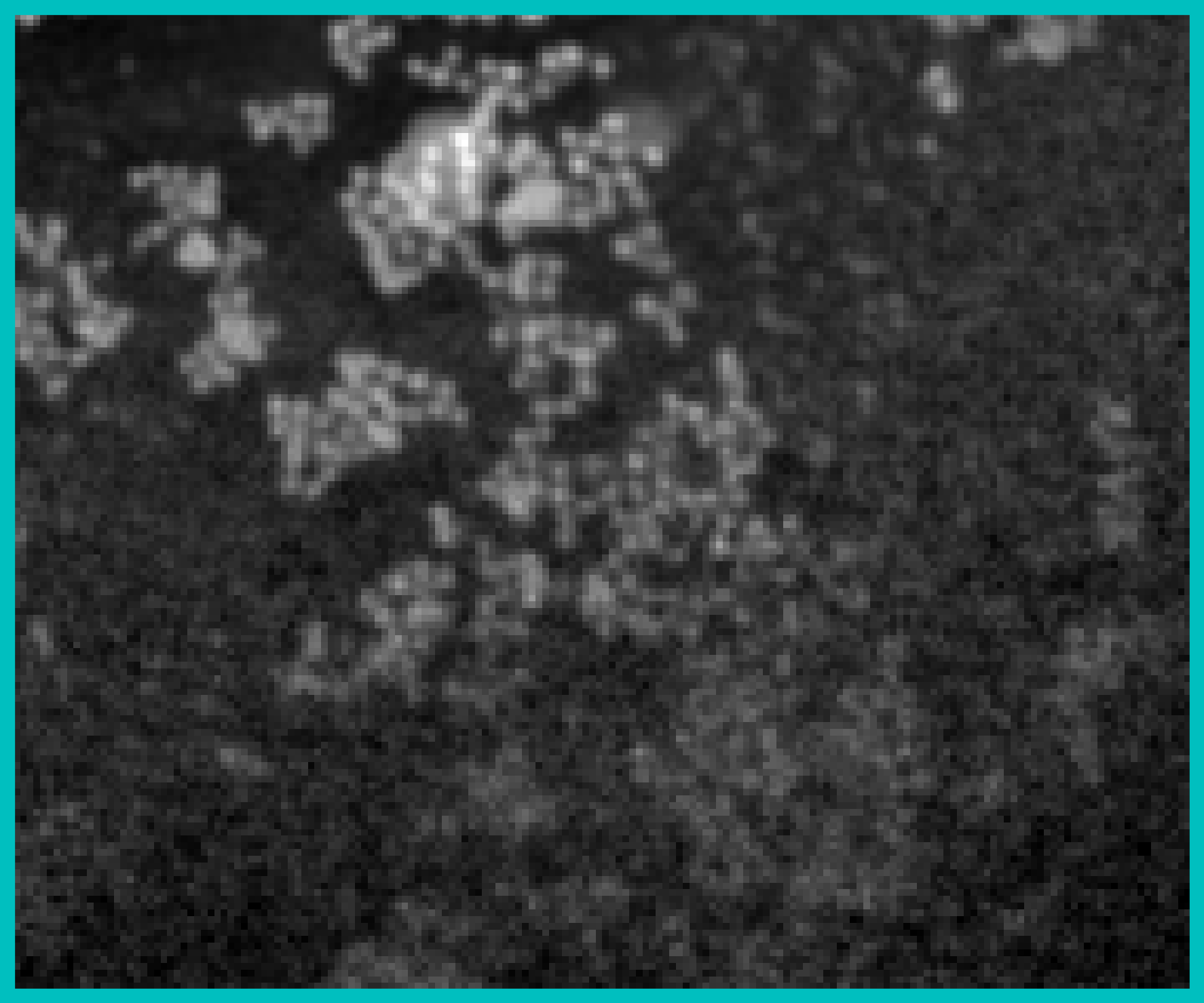}}
  \vspace{0.1cm}
  \centerline{(d) PMACE}\medskip  
\end{minipage}
\hfill
\hfill
\vspace{-0.1cm}
\caption{\textbf{Phases} (in radians) of complex reconstructions from measured data. The top row shows the full field of view along with outlines of insets displayed in rows 2-4.  All methods reconstruct fine detail in the central region (row 3, outlined in yellow).  ePIE and AWF yield artifacts near the edge of the field of view, while SHARP and PMACE produce more uniform results throughout the field of view.}
\label{fig:GoldBalls_phase}
\vspace{-0.8cm}
\end{figure*}

\begin{figure*}[!b]
\begin{minipage}[b]{.21\linewidth}
  \centering
  \centerline{\includegraphics[height=4cm]{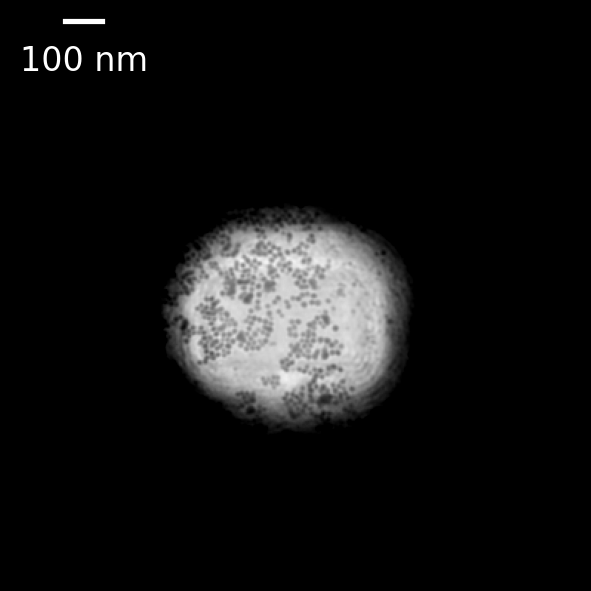}}
  \vspace{0.1cm}
  \centerline{(a) ePIE}\medskip
\end{minipage}
\hfill
\begin{minipage}[b]{.21\linewidth}
  \centering
  \centerline{\includegraphics[height=4cm]{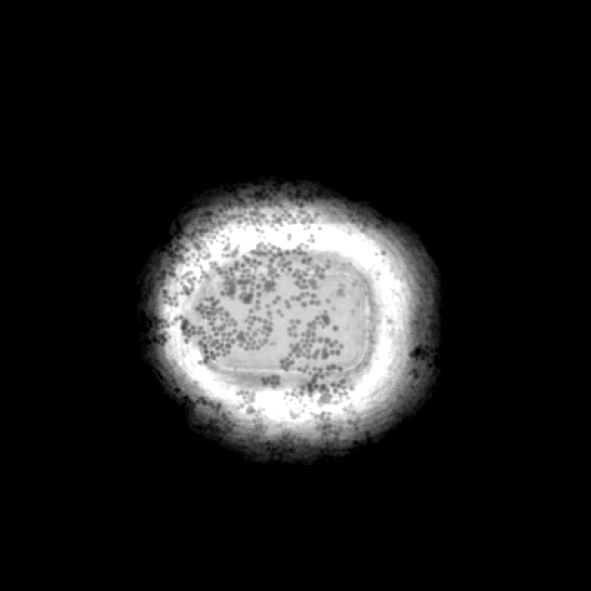}}
  \vspace{0.1cm}
  \centerline{(b) AWF}\medskip
\end{minipage}
\hfill
\begin{minipage}[b]{.21\linewidth}
  \centering
  \centerline{\includegraphics[height=4cm]{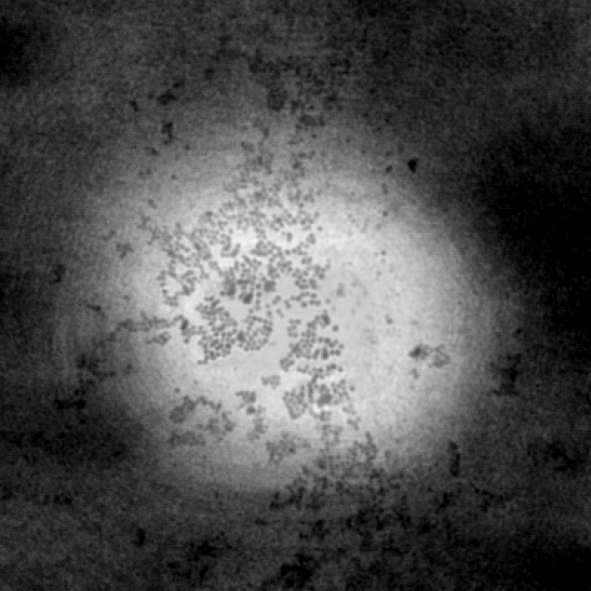}}
  \vspace{0.1cm}
  \centerline{(c) SHARP}\medskip  
\end{minipage}
\hfill
\begin{minipage}[b]{.21\linewidth}
  \centering
  \centerline{\includegraphics[height=4cm]{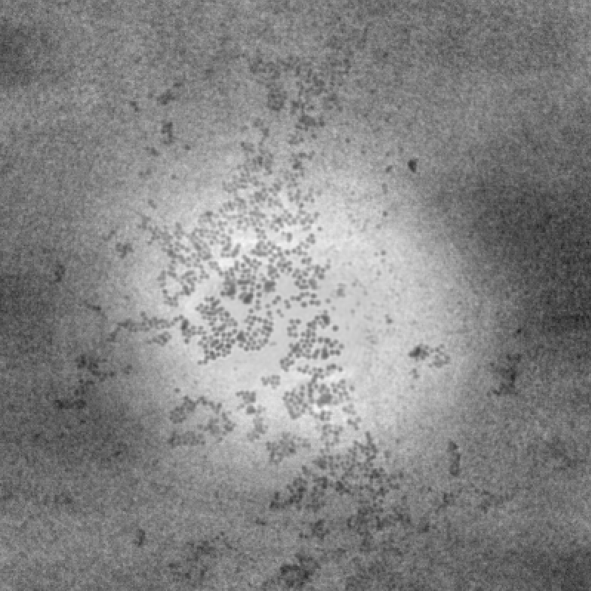}}
  \vspace{0.1cm}
  \centerline{(d) PMACE}\medskip  
\end{minipage}
\hfill
\hspace{0.06cm}
\begin{minipage}[b]{.02\linewidth}
  \centering
  \centerline{\includegraphics[height=4.1cm]{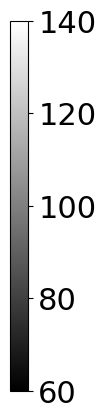}}
  \centerline{}\medskip
\end{minipage}
\hfill
\vspace{-0.1cm}
\caption{\textbf{Magnitudes} of complex reconstructions from measured data. PMACE captures fine detail in the largest field of view relative to the other methods.}
\vspace{0.5cm}
\label{fig:GoldBalls_mag}
\end{figure*}

The reconstructed phase and magnitude images are shown in Figures~\ref{fig:GoldBalls_phase} and~\ref{fig:GoldBalls_mag}.
We cropped the reconstructions to show the center $400 \times 400$ region of the image that was fully measured by the probes, and we used color outlines to highlight zoomed-in regions of the phase images.
As can be observed, ePIE and AWF produce good image quality in the central region (yellow), but their image quality degrades in the red and blue regions.
SHARP and PMACE produce good results in all three regions, with PMACE reconstructing the gold balls in the largest field of view.
PMACE also appears to have slightly higher sharpness/contrast than SHARP, but since we do not have ground truth, it is difficult to know if the detail in the background is noise or valid detail.

\begin{table}[!htbp]
    \centering
    \begin{tabular}{|c|c|c|c|c|} \hline 
         Algorithm & ePIE & AWF & SHARP & PMACE \\ \hline 
         NRMSE & 0.1102 & 0.1041 & 0.0914 & {\bf 0.0912} \\ \hline 
    \end{tabular}
    \caption{Forward Propagated NRMSE for each of the reconstruction algorithms. Notice that PMACE results in the lowest forward propagated error.}
    \label{tab:ForwardPropagatedNRMSE}
\end{table}

Table~\ref{tab:ForwardPropagatedNRMSE} tabulates a measure of reconstruction quality for the ePIE, AWF, SHARP and PMACE reconstructions on the gold balls data set.
Since there is no ground truth, we forward propagated the reconstructed complex images to the detector plane following \eqref{eq: mag_measurement} and calculate the NRMSE in the magnitude relative to the square-rooted data. 
In principle, a lower NRMSE at the detector plane implies a better fit to measurements; however, one must be careful since this metric does not account for overfitting of data.
Notice that the PMACE reconstruction provides the best fit to the measured data.   
Hence these results support the claim that PMACE provides state-of-the-art quality under high-overlap conditions for measured data.

\section{Conclusion}
In this paper, we introduce a theoretical foundation for the PMACE algorithm and show how it can be applied to the specific problem of ptychographic reconstruction.
The PMACE algorithm is based on a set of agents each of which only acts on a small portion of the full reconstruction.
This makes PMACE more suitable for large distributed inverse problems in which measured data only acts on a portion of the total image.
As with MACE, the PMACE algorithm solves an equilibrium problem, which does not always correspond to an optimization problem.
This makes PMACE more flexible for use in applications than traditional regularized inversion.

We applied PMACE to ptychography on both synthetic and measured data and found that it had improved image quality relative to competing algorithms, particularly for the case of large spacing between probe positions. 
This last point implies that our method requires less data for accurate reconstruction and can thus reduce acquisition time for a given level of reconstruction quality.  

Finally, we note that the PMACE algorithm can be applied to other inverse problems in which each measurement only depends on a portion of the reconstruction.
For example, this may be appropriate in tomography applications, such as laminography, in which each view corresponds to a portion of a larger planner object.

{\appendices
{
\section{Proofs}

\subsection{Proof of Theorem~\ref{thm: using_prox_equal_optimization_problem}}
\extraproof{
\begin{theorem}
Let $F_j, j = 0, \cdots, J-1$ denote the proximal map function of a differentiable convex function $f_j$ as specified in \eqref{eq: proximal_map_function}.  Then $\bv^*$ is a solution to the PMACE equation in \eqref{eq:pmace-eqns} if and only if $x^* = \bar{x}(\bv^*)$ satisfies
\begin{equation}
    \sum_{j=0}^{J-1} P_j^T W \nabla f_j(P_j x^*) = 0.
\end{equation}
\end{theorem}
}
\begin{proof}
The assumptions on $f_j$ and $F_j$ imply~\cite[Sec. 6.1]{boyd2016primer} that $F_j(v_j) = (I + \sigma^2 \nabla f_j)^{-1}(v_j)$ and hence $v_j = (I + \sigma^2 \nabla f_j)(F_j(v_j))$.  With this in mind, suppose $\bv^*$ satisfies \eqref{eq:pmace-eqns}.  By definition of $\bF$, $\bG^P$, and $x^*$, this means
\begin{equation} \label{eq:FjPj}
    F_j(v_j^*) = P_j \bar{x}(\bv^*) = P_j x^*
\end{equation}
for all $j$.  Using $F_j = (I + \sigma^2 \nabla f_j)^{-1}$ and then applying $(I + \sigma^2 \nabla f_j)$ to both sides gives
\begin{equation} \label{eq:vjPj}
    v_j^* = P_j x^* + \sigma^2 \nabla f_j(P_j x^*) \ .
\end{equation}
Using this with the definition of $x^*$ gives
\begin{equation} 
    x^* = \bar{x}(\bv^*) = \Lambda^{-1} \sum_j P_j^T W(P_j x^* + \sigma^2 \nabla f_j(P_j x^*)) \ ,
\end{equation}
so
\begin{equation} \label{eq:Lambda-x-star}
    \Lambda x^* =  \sum_j P_j^T W P_j x^* + \sigma^2 \sum_j P_j^T W \nabla f_j(P_j x^*)) \ .
\end{equation}
Since the first term on the right hand side is $\Lambda x^*$, this simplifies to 
\begin{equation} \label{eq:sigma-vec-field}
    \sigma^2 \sum_j P_j^T W \nabla f_j(P_j x^*)) = 0 \ ,
\end{equation}
which gives \eqref{eq:vec-field}.

For the converse, if $x^*$ satisfies \eqref{eq:vec-field}, then by taking \eqref{eq:vjPj} as the definition of $v_j^*$, each step above is reversible to give \eqref{eq:FjPj} and hence \eqref{eq:pmace-eqns}.
\end{proof}

\subsection{Proof of Theorem~\ref{thm: using_prox_equil_problem}}
\extraproof{
\begin{theorem}
Let $f_j$ and $F_j$ be as in Theorem~\ref{thm: using_prox_equal_optimization_problem}.  
\begin{enumerate}
    \item If $W = rI$ for scalar $r$.  Then PMACE is equivalent to optimization in that $\bv^*$ is a solution to the PMACE equation in \eqref{eq:pmace-eqns} if and only if $x^* = \bar{x}(\bv^*)$ satisfies
    \begin{equation}
    x^* = \arg\min_{x}\left \{\sum_{j=0}^{J-1}f_{j}(P_j x) \right \} \ .
    \end{equation}
    \item For generic diagonal $W$ and $f_j$, the PMACE formulation does not naturally arise as an optimization problem.  That is, the vector field in \eqref{eq:vec-field} is not a conservative vector field, and hence is not the gradient field of a potential function. 
\end{enumerate}
\end{theorem}
}
\begin{proof}
1.  Note that $\bG^P$ is independent of $r>0$ when $W = rI$, so we may assume $W = I$.  A straightforward modification to the proof of~\cite[Theorem 2]{buzzard2018plug}, shows that if $\bF$ is defined using the $F_j$, then a PMACE solution $\bv^*$ of \eqref{eq:pmace-eqns} is equivalent using $x^* = \bar{x}(\bv^*)$, $\lambda^*_j = v_j^* - P_j x^*$ to solutions  $(x^{*}, \lambda_{j}^*)$ that solve the following equations 
\begin{equation}
\label{eq: pmace_equilibrium_eq_data_fidelity}
    F_j(P_j x^* + \lambda_{j}^*) = P_j x^*
\end{equation}
\begin{equation}
\label{eq: pmace_equilibrium_condition}
     \sum_{j} P_{j}^{T} \lambda^{*}_{j} = 0 \ .
\end{equation}
These conditions are equivalent to the first-order optimality conditions for the constrained optimization problem 
\begin{equation} \label{eq:x-split}
x^* = \arg\min_{x, u_j}\left \{\sum_{j=0}^{J-1}f_{j}(u_j) \right \} \ \text{s.t.} \ u_j = P_j x \text{ for all } j \ ,
\end{equation}
which is equivalent to \eqref{eq: ml_est}. 

2. As noted above, the conservative property is equivalent to the symmetry of the Hessian matrix in \eqref{eq:grad-vec-field}, which is equivalent to 
\begin{equation} \label{eq:Hessian}
\sum_j P_j^T(W H_j(x) - H_j(x) W) P_j = 0 \ ,
\end{equation}
for all $x$, where $H_j(x)$ is $H f_j$ evaluated at $P_j x$.  After an arbitrarily small perturbation of the diagonal elements of $W$, we may assume that these diagonal entries are all distinct and positive. Then, if the left hand side of \eqref{eq:Hessian} is not identically 0, we are done.  Otherwise, by an $\epsilon$-small perturbation of $f_0$ (in the sup-norm on a ball large enough to include all images of interest), we can replace $H_0$ in \eqref{eq:Hessian} by $H_0 + \epsilon H$, with at least two off-diagonal entries of $H$ nonzero.  The assumption that $P_j P_j^T = I$ implies that $P_0^T$ is injective, so with $H_0 + \epsilon H$ in place of $H_0$, it is sufficient to show that $WH \neq HW$.   Since $W$ is diagonal with distinct, positive diagonal entries and $H$ is symmetric, they commute if and only if $H$ is also diagonal, which contradicts the choice of $H$.  Hence after an arbitrarily small perturbation of $W$ and/or $f_0$, there is no potential function associated with the vector field \eqref{eq:vec-field}.  
\end{proof}

\subsection{Proof of Theorem~\ref{thm: pmace_sol_equal_v_goes_zeros}}
\extraproof{
\begin{theorem} 
The solution of~\eqref{eq:pmace-eqns} is equivalent to solving the inclusion problem $0 \in V(x^*)$, where $V$ is the set-valued vector field 
\begin{equation}  %
   V(x) = \sum_{j=0}^{J-1} P_j^T W \left( P_j x  - F_j^{-1}(P_j x)\right) \ . 
\end{equation}
\end{theorem}
}
\begin{proof}
Let $\bv^*$ denote the solution of~\eqref{eq:pmace-eqns} and $x^* = \bar{x}(\bv^*)$. Since $\bG^P(\bx^*) = \bG^P( \bG^P(\bx^*) )$, $\bv^*$ satisfies \eqref{eq:pmace-eqns} if and only if $F_j(v_j^*) = P_j x^*$ for all $j$ and
\begin{equation}
    \bG^P(\bv^*) = \bG^P( \bF(\bv^*) ) \ .
\end{equation}
This is equivalent to $v_j \in F_j^{-1}(P_j x^*)$ for all $j$ and 
\begin{equation}
    P_{i} \Lambda^{-1} \sum_{j=0}^{J-1} P_{j}^{T} W v_{j}^{*} = P_{i} \Lambda^{-1} \sum_{j=0}^{J-1} P_{j}^{T} W F_{j} (v_j^{*}) \ ,
\end{equation}
 for all $i$.  
Since the matrix obtained by stacking the $P_i$ is injective, this is equivalent to $v_j \in F_j^{-1}(P_j x^*)$ for all $j$ and 
\begin{equation}
    \sum_{j=0}^{J-1} P_j^T W \left( v_j^*  - F_j(v_j^*)\right) = 0 \ .
\end{equation}
Replacing both $v_j^*$ with $F_j^{-1}(P_j x^*)$ gives the theorem.  
\end{proof}
}

\subsection{Proof that the PMACE forward operator is invertible}
\begin{proof}
Recall that the PMACE forward operator from \eqref{eq:Fjv} is
\begin{align}
    F_{j} (v_{j})
        & = (1 - \alpha ) v_{j} + \alpha D^{-1}\mathcal{F}^{*} \left ( y_{j} \frac{\mathcal{F} D v_{j}}{|\mathcal{F} D v_{j}|} \right ) \ .
\end{align}
We first multiply both sides of \eqref{eq:Fjv} by $\mathcal{F} D$ to get
\begin{equation} \label{eq:FDFj}
\mathcal{F} D F_j(v_j) = \left ( \left ( 1-\alpha \right ) |\mathcal{F} D v_{j}| + \alpha y_j \right ) \frac{\mathcal{F} D v_{j}}{|\mathcal{F} D v_{j}|} \ .
\end{equation}
Taking absolute value of both sides and using $y_j \geq 0$ allows us to solve for $|\mathcal{F} D v_{j}|$ in terms of $F_j(v_j)$, $y_j$, and $\alpha$. 
This leaves $\mathcal{F} D v_{j}$ as the only unknown in \eqref{eq:FDFj}, so we solve for that and then multiply by $D^{-1} \mathcal{F}^*$ to get $v_j$.
\end{proof}
}

\bibliographystyle{IEEEtran}
\bibliography{references}

\begin{thebibliography}{10}
\providecommand{\url}[1]{#1}
\csname url@samestyle\endcsname
\providecommand{\newblock}{\relax}
\providecommand{\bibinfo}[2]{#2}
\providecommand{\BIBentrySTDinterwordspacing}{\spaceskip=0pt\relax}
\providecommand{\BIBentryALTinterwordstretchfactor}{4}
\providecommand{\BIBentryALTinterwordspacing}{\spaceskip=\fontdimen2\font plus
\BIBentryALTinterwordstretchfactor\fontdimen3\font minus \fontdimen4\font\relax}
\providecommand{\BIBforeignlanguage}[2]{{%
\expandafter\ifx\csname l@#1\endcsname\relax
\typeout{** WARNING: IEEEtran.bst: No hyphenation pattern has been}%
\typeout{** loaded for the language `#1'. Using the pattern for}%
\typeout{** the default language instead.}%
\else
\language=\csname l@#1\endcsname
\fi
#2}}
\providecommand{\BIBdecl}{\relax}
\BIBdecl

\bibitem{dabov2007}
K.~Dabov, A.~Foi, V.~Katkovnik, and K.~Egiazarian, ``{Image Denoising by Sparse 3-D Transform-Domain Collaborative Filtering},'' \emph{IEEE Trans. Image Process.}, vol.~16, no.~8, pp. 2080--2095, 2007.

\bibitem{zhang2017beyond}
K.~Zhang, W.~Zuo, Y.~Chen, D.~Meng, and L.~Zhang, ``Beyond a {G}aussian denoiser: Residual learning of deep {CNN} for image denoising,'' \emph{IEEE transactions on image processing}, vol.~26, no.~7, pp. 3142--3155, 2017.

\bibitem{buzzard2018plug}
G.~T. Buzzard, S.~H. Chan, S.~Sreehari, and C.~A. Bouman, ``Plug-and-play unplugged: Optimization-free reconstruction using consensus equilibrium,'' \emph{SIAM Journal on Imaging Sciences}, vol.~11, no.~3, pp. 2001--2020, 2018.

\bibitem{venkatakrishnan_plug-and-play_2013}
S.~V. Venkatakrishnan, C.~A. Bouman, and B.~Wohlberg, ``Plug-and-{Play} priors for model based reconstruction,'' in \emph{2013 {IEEE} {Global} {Conference} on {Signal} and {Information} {Processing}}, Dec. 2013, pp. 945--948.

\bibitem{sreehari2016plug}
S.~Sreehari, S.~V. Venkatakrishnan, B.~Wohlberg, G.~T. Buzzard, L.~F. Drummy, J.~P. Simmons, and C.~A. Bouman, ``Plug-and-play priors for bright field electron tomography and sparse interpolation,'' \emph{Transactions on Computational Imaging}, vol.~2, pp. 408--423, 2016.

\bibitem{bouman2022foundations}
\BIBentryALTinterwordspacing
C.~A. Bouman, \emph{Foundations of Computational Imaging: A Model-Based Approach}.\hskip 1em plus 0.5em minus 0.4em\relax Philadelphia, PA: Society for Industrial and Applied Mathematics, 2022. [Online]. Available: \url{https://epubs.siam.org/doi/abs/10.1137/1.9781611977134}
\BIBentrySTDinterwordspacing

\bibitem{rodenburg2008ptychography}
J.~M. Rodenburg, ``Ptychography and related diffractive imaging methods,'' \emph{Advances in imaging and electron physics}, vol. 150, pp. 87--184, 2008.

\bibitem{rodenburg2019ptychography}
J.~Rodenburg and A.~Maiden, ``Ptychography,'' \emph{Springer Handbook of Microscopy}, pp. 819--904, 2019.

\bibitem{wilke2012hard}
R.~Wilke, M.~Priebe, M.~Bartels, K.~Giewekemeyer, A.~Diaz, P.~Karvinen, and T.~Salditt, ``Hard x-ray imaging of bacterial cells: nano-diffraction and ptychographic reconstruction,'' \emph{Optics express}, vol.~20, no.~17, pp. 19\,232--19\,254, 2012.

\bibitem{trtik2013density}
P.~Trtik, A.~Diaz, M.~Guizar-Sicairos, A.~Menzel, and O.~Bunk, ``Density mapping of hardened cement paste using ptychographic x-ray computed tomography,'' \emph{Cement and Concrete Composites}, vol.~36, pp. 71--77, 2013.

\bibitem{guizar2014high}
M.~Guizar-Sicairos, I.~Johnson, A.~Diaz, M.~Holler, P.~Karvinen, H.-C. Stadler, R.~Dinapoli, O.~Bunk, and A.~Menzel, ``High-throughput ptychography using eiger: scanning x-ray nano-imaging of extended regions,'' \emph{Optics express}, vol.~22, no.~12, pp. 14\,859--14\,870, 2014.

\bibitem{shapiro2017ptychographic}
D.~A. Shapiro, R.~Celestre, P.~Denes, M.~Farmand, J.~Joseph, A.~Kilcoyne, S.~Marchesini, H.~Padmore, S.~V. Venkatakrishnan, T.~Warwick \emph{et~al.}, ``Ptychographic imaging of nano-materials at the advanced light source with the nanosurveyor instrument,'' in \emph{Journal of Physics: Conference Series}, vol. 849, no.~1, 2017, p. 012028.

\bibitem{fienup1982phase}
J.~R. Fienup, ``Phase retrieval algorithms: a comparison,'' \emph{Applied Optics}, vol.~21, no.~15, pp. 2758--2769, 1982.

\bibitem{rodenburg2004phase}
J.~M. Rodenburg and H.~M. Faulkner, ``A phase retrieval algorithm for shifting illumination,'' \emph{Applied Physics Letters}, vol.~85, no.~20, pp. 4795--4797, 2004.

\bibitem{maiden2009improved}
A.~M. Maiden and J.~M. Rodenburg, ``An improved ptychographical phase retrieval algorithm for diffractive imaging,'' \emph{Ultramicroscopy}, vol. 109, no.~10, pp. 1256--1262, 2009.

\bibitem{thibault_reconstructing_2013}
\BIBentryALTinterwordspacing
P.~Thibault and A.~Menzel, ``Reconstructing state mixtures from diffraction measurements,'' \emph{Nature}, vol. 494, no. 7435, pp. 68--71, Feb. 2013. [Online]. Available: \url{https://doi.org/10.1038/nature11806}
\BIBentrySTDinterwordspacing

\bibitem{maiden2017further}
A.~Maiden, D.~Johnson, and P.~Li, ``Further improvements to the ptychographical iterative engine,'' \emph{Optica}, vol.~4, no.~7, pp. 736--745, 2017.

\bibitem{bunk2008influence}
O.~Bunk, M.~Dierolf, S.~Kynde, I.~Johnson, O.~Marti, and F.~Pfeiffer, ``Influence of the overlap parameter on the convergence of the ptychographical iterative engine,'' \emph{Ultramicroscopy}, vol. 108, no.~5, pp. 481--487, 2008.

\bibitem{candes2015phase}
E.~J. Candes, X.~Li, and M.~Soltanolkotabi, ``Phase retrieval via {W}irtinger flow: Theory and algorithms,'' \emph{IEEE Transactions on Information Theory}, vol.~61, no.~4, pp. 1985--2007, 2015.

\bibitem{xu2018accelerated}
R.~Xu, M.~Soltanolkotabi, J.~P. Haldar, W.~Unglaub, J.~Zusman, A.~F. Levi, and R.~M. Leahy, ``Accelerated {W}irtinger flow: A fast algorithm for ptychography,'' \emph{arXiv preprint arXiv:1806.05546}, 2018.

\bibitem{odstrvcil2018iterative}
M.~Odstr{\v{c}}il, A.~Menzel, and M.~Guizar-Sicairos, ``Iterative least-squares solver for generalized maximum-likelihood ptychography,'' \emph{Optics Express}, vol.~26, no.~3, pp. 3108--3123, 2018.

\bibitem{soulez2016proximity}
F.~Soulez, {\'E}.~Thi{\'e}baut, A.~Schutz, A.~Ferrari, F.~Courbin, and M.~Unser, ``Proximity operators for phase retrieval,'' \emph{Applied optics}, vol.~55, no.~26, pp. 7412--7421, 2016.

\bibitem{yan2020ptychographic}
H.~Yan, ``Ptychographic phase retrieval by proximal algorithms,'' \emph{New Journal of Physics}, vol.~22, no.~2, p. 023035, 2020.

\bibitem{marchesini2016alternating}
S.~Marchesini, Y.-C. Tu, and H.-t. Wu, ``Alternating projection, ptychographic imaging and phase synchronization,'' \emph{Applied and Computational Harmonic Analysis}, vol.~41, no.~3, pp. 815--851, 2016.

\bibitem{marchesini2016sharp}
S.~Marchesini, H.~Krishnan, B.~J. Daurer, D.~A. Shapiro, T.~Perciano, J.~A. Sethian, and F.~R. Maia, ``{SHARP}: a distributed {GPU}-based ptychographic solver,'' \emph{Journal of Applied Crystallography}, vol.~49, no.~4, pp. 1245--1252, 2016.

\bibitem{luke2004relaxed}
D.~R. Luke, ``Relaxed averaged alternating reflections for diffraction imaging,'' \emph{Inverse Problems}, vol.~21, no.~1, p.~37, 2004.

\bibitem{sridhar2020distributed}
V.~Sridhar, X.~Wang, G.~T. Buzzard, and C.~A. Bouman, ``Distributed iterative ct reconstruction using multi-agent consensus equilibrium,'' \emph{IEEE Transactions on Computational Imaging}, vol.~6, pp. 1153--1166, 2020.

\bibitem{williamson2004multivariable}
R.~E. Williamson and H.~F. Trotter, ``Multivariable mathematics: Linear algebra, calculus,'' \emph{Differential Equations}, vol.~2, no.~4, 2004.

\bibitem{combettes_monotone_2018}
\BIBentryALTinterwordspacing
P.~L. Combettes, ``Monotone operator theory in convex optimization,'' \emph{Mathematical Programming}, vol. 170, no.~1, pp. 177--206, Jul. 2018. [Online]. Available: \url{https://doi.org/10.1007/s10107-018-1303-3}
\BIBentrySTDinterwordspacing

\bibitem{Bauschke2011ConvexAA}
H.~H. Bauschke and P.~L. Combettes, ``Convex analysis and monotone operator theory in {H}ilbert spaces,'' in \emph{CMS Books in Mathematics}, 2011.

\bibitem{venkatakrishnan2013plug}
S.~V. Venkatakrishnan, C.~A. Bouman, and B.~Wohlberg, ``Plug-and-play priors for model based reconstruction,'' in \emph{IEEE Global Conference on Signal and Information Processing}.\hskip 1em plus 0.5em minus 0.4em\relax IEEE, 2013, pp. 945--948.

\bibitem{godard2012noise}
P.~Godard, M.~Allain, V.~Chamard, and J.~Rodenburg, ``Noise models for low counting rate coherent diffraction imaging,'' \emph{Optics Express}, vol.~20, no.~23, pp. 25\,914--25\,934, 2012.

\bibitem{huang2014optimization}
X.~Huang, H.~Yan, R.~Harder, Y.~Hwu, I.~K. Robinson, and Y.~S. Chu, ``Optimization of overlap uniformness for ptychography,'' \emph{Optics Express}, vol.~22, no.~10, pp. 12\,634--12\,644, 2014.

\bibitem{zhou2020low}
L.~Zhou, J.~Song, J.~S. Kim, X.~Pei, C.~Huang, M.~Boyce, L.~Mendon{\c{c}}a, D.~Clare, A.~Siebert, C.~S. Allen \emph{et~al.}, ``Low-dose phase retrieval of biological specimens using cryo-electron ptychography,'' \emph{Nature Communications}, vol.~11, no.~1, pp. 1--9, 2020.

\bibitem{marchesini2017ptychography}
S.~Marchesini, ``Ptychography gold ball example dataset,'' Lawrence Berkeley National Lab. (LBNL), Berkeley, CA, USA, Tech. Rep., 2017.

\bibitem{boyd2016primer}
E.~K. Ryu and S.~Boyd, ``A primer on monotone operator methods survey,'' \emph{Appl. Comput. Math.}, vol.~15, no.~1, pp. 3--43, 2016.

\end{thebibliography}

\newpage

\vfill

\end{document}